\documentclass[journal]{IEEEtran}

\usepackage{hyperref}
\usepackage{enumerate}
\usepackage{amsmath,amssymb,amsthm}

\newtheorem{theorem}{Theorem}
\newtheorem{prop}{Proposition}
\newtheorem{lemma}{Lemma}
\newtheorem{corollary}{Corollary}
\newtheorem{remark}{Remark}
\newtheorem{assumption}{Assumption}

\newcommand{\be}{\begin{equation}}
\newcommand{\ee}{\end{equation}}
\newcommand{\ba}{\begin{array}}
\newcommand{\ea}{\end{array}}
\newcommand{\bea}{\begin{eqnarray}}
\newcommand{\eea}{\end{eqnarray}}

\newcommand{\diag}{{\mbox{diag}}}

\newcommand{\tran}{^{\mbox{\scriptsize T}}}  

\newcommand{\vbar}{\raisebox{.17ex}{\rule{.04em}{1.35ex}}}
\newcommand{\vbarind}{\raisebox{.01ex}{\rule{.04em}{1.1ex}}}
\newcommand{\D}{\ifmmode {\rm I}\hspace{-.2em}{\rm D} \else ${\rm I}\hspace{-.2em}{\rm D}$ \fi}
\newcommand{\T}{\ifmmode {\rm I}\hspace{-.2em}{\rm T} \else ${\rm I}\hspace{-.2em}{\rm T}$ \fi}
\newcommand{\B}{\ifmmode {\rm I}\hspace{-.2em}{\rm B} \else \mbox{${\rm I}\hspace{-.2em}{\rm B}$} \fi}
\newcommand{\Hil}{\ifmmode {\rm I}\hspace{-.2em}{\rm H} \else \mbox{${\rm I}\hspace{-.2em}{\rm H}$} \fi}
\newcommand{\C}{\ifmmode \hspace{.2em}\vbar\hspace{-.31em}{\rm C} \else \mbox{$\hspace{.2em}\vbar\hspace{-.31em}{\rm C}$} \fi}
\newcommand{\Cind}{\ifmmode \hspace{.2em}\vbarind\hspace{-.25em}{\rm C} \else \mbox{$\hspace{.2em}\vbarind\hspace{-.25em}{\rm C}$} \fi}
\newcommand{\Q}{\ifmmode \hspace{.2em}\vbar\hspace{-.31em}{\rm Q} \else \mbox{$\hspace{.2em}\vbar\hspace{-.31em}{\rm Q}$} \fi}
\newcommand{\Z}{\ifmmode {\rm Z}\hspace{-.28em}{\rm Z} \else ${\rm Z}\hspace{-.38em}{\rm Z}$ \fi}

\renewcommand{\vec}[1]{{\bf{#1}}}     

\usepackage{graphicx}
\usepackage{epstopdf}
\usepackage{caption}
\usepackage{subcaption}
\newcommand{\R}{\mathbb{R}}
\newcommand{\N}{\mathbb{N}}
\newcommand{\Blambda}{\boldsymbol\lambda}
\newcommand{\Bmu}{\boldsymbol\mu}

\usepackage{color,soul}
\definecolor{lightblue}{rgb}{.90,.95,1}
\sethlcolor{lightblue}
\usepackage{tikz}
\usetikzlibrary{shapes,arrows}

\newcommand{\bec}[1]{\bar{\vec{#1}}}
\newcommand{\hec}[1]{\hat{\vec{#1}}}

\newcommand{\dist}{\texttt{dist}}
\newcommand{\fes}{\texttt{fes}}
\newcommand{\sign}{\texttt{sign}}

\title{Voltage Control Using Limited Communication}
\author{Sindri  Magn\'{u}sson, Guannan Qu, Carlo Fischione, and  Na Li
\thanks{This work was supported by the VR Chromos Project, NSF EPCN 1608509, NSF CAREER 1553407, and ARPA-E NODES.}
\thanks{S.  Magn\'{u}sson, G. Qu, and N. Li is with the School of Engineering and Applied Sciences, Harvard University, Cambridge, MA, USA (Email: {sindrim@seas.harvard.edu} {gqu@g.harvard.edu},{nali@seas.harvard.edu}).
}
\thanks{C. Fischione are with the School of Electrical Engineering and Computer Science, KTH Royal Institute of Technology, Stockholm, Sweden (Email:{carlofi@kth.se}).} 
}

\begin{document}


%
%
%


\maketitle

\begin{abstract}

 In electricity  distribution networks, the increasing penetration of renewable energy generation necessitates faster and more sophisticated voltage controls.  Unfortunately, recent research shows that local voltage control fails in achieving the desired regulation, unless there is communication between the controllers. However, the communication infrastructure for distribution systems is less reliable and less ubiquitous as compared to that for the bulk transmission system. In this paper, we design distributed voltage control that uses limited communication. That is, only neighboring buses need to communicate a few bits between each other for each control step. We investigate how these controllers can achieve the desired asymptotic behavior of the voltage regulation and  we provide upper bounds on the number of bits that are needed to ensure a predefined accuracy of the regulation. Finally, we illustrate the results by numerical simulations. 
     
\end{abstract}




 \section{Introduction}
   
 There is an increasing penetration of distributed energy resources such as 
 renewable energy in distribution networks. 
  Unfortunately, such a penetration causes  faster voltage fluctuations than what
   today's distribution networks can handle, see~\cite{carvalho2008distributed}.  
Therefore, to avoid overloading the distribution networks, the integration of renewable energy resources must be accompanied by faster and more sophisticated voltage regulation.

 These challenges have motivated a growing research interest in voltage control, where fast voltage fluctuations are regulated through real-time reactive power injections to ensure that the voltage is maintained within an acceptable range. Such fast voltage control can be implemented in the emerging power devices such as inverters. The research efforts have focused on two main directions: local and  distributed control strategies.  In the local voltage control, control devices at each bus update the reactive power injections using only locally available information, such as local voltage measurements, 
see \cite{farivar2013equilibrium,Li2014,Zhu_2015} and references therein.  On the other hand, in distributed voltage control schemes, control devices at each bus determine the reactive power injection with additional information communicated from its neighboring buses in the distribution network,  see \cite{Zhang_2015,vsulc2014optimal,Bolognani2013,Bolognani2015,liu2017hybrid}. 
 Local control strategies have the obvious advantage over distributed ones in that they do not rely  on communication.  
 However, even though  local control strategies perform well in some cases, they may fail to ensure that the voltage is maintained within the accepted range in some cases, as proved by the impossibility result in~\cite{Cavraro2016}. 
 Therefore, communication among the local controllers is needed to guarantee the performance of voltage regulation. 

 The communication capabilities of today's distribution networks generally suffer from low data rates,~\cite{Yan_2013,Galli2011grid}.  
 To compensate for this deficiency, power system operators and industries are currently investing heavily in integrating the distribution networks with a sophisticated communication infrastructure. 
 However, even with the promising  capabilities of the future  low latency  networks, fast real-time control applications, like voltage control, rely on short packages that carry coarsely quantized information,~\cite{Durisi_2016}. 
 Therefore, it is important to develop voltage control with very limited communication    
   for early integration of renewable resources using today's grid limited communication capabilities and  also for sustainable developments of the future smart grid.

   Networked control systems with limited bandwidth communication have been well investigated~\cite{brockett2000quantized,nair2007feedback,zhang2013network}.
   However,   standard algorithms used for stabilization of networked control systems are generally not applicable to the voltage regulation problem. 
   Instead, most distributed voltage control algorithms in the literature are based on the tools of distributed optimization, e.g.~\cite{Zhang_2015,vsulc2014optimal,Bolognani2013,Bolognani2015,liu2017hybrid}. 
  Even though distributed optimization with limited bandwidth has received some attention~\cite{Rabbat_2005,nedic2008distributed,Magnusson_2016,Magnusson2018}, none of these works fit the characteristics of the voltage control problem.

%
%
   
 {\color{black}   
   Other ways of limiting the communication of voltage control algorithms have been considered in the literature. 
   The work in~\cite{olivier2016active} proposes that the systems only communicates when there is an overvoltage in the system. 
    Whereas the method in~\cite{fan2017distributed} proposes that buses communicate their physical state only whent he difference between their physical state and previously communicated state is too large.
      Unlike our work, these papers do not consider bandwidth limited communication.
%
}



 In this paper, we study a distributed voltage control where only a few bits of communication between neighboring buses are needed during each control step.  
 In particular,  the voltage control device on each bus determines the reactive power injection based on its local voltage measurement and current reactive power injection, in addition to a few bits of information communicated from its physical neighbors.
  We show that the algorithm can regulate the voltages to an acceptable range, for any  predefined accuracy, in a finite number of iterations. 
    We also provide an upper bound  on the number of communicated bits  (in the worst case) that are needed to ensure a predefined accuracy of the desired voltage level. Though the theoretical analysis is based on a linearized power flow model that is applicable for radial distribution network developed in \cite{Baran1989}, we use the nonlinear power flow model to numerically test the algorithm. Moreover, we test the performance under both static and dynamic operating conditions. Both the theoretical and numerical results confirm that reducing the communication to a few bits do not sacrifice the control performance. This opens up a large flexibility for implementing voltage control in practice. For instance it enables the use of power lines as the communication media for electricity distribution systems despite the limited bandwidth of power line communication~\cite{Galli2011grid}.  
%
%
%
%
%

    Preliminary studies of this work  appeared in~\cite{SindriIFAC2017}. 
   However,  in~\cite{SindriIFAC2017} most of the proofs are omitted (but they appear here) and all the numerical results especially with the nonlinear power flow model here are new. 
   Moreover, the presentation and theoretical results here have been largely improved compared to~\cite{SindriIFAC2017}. 
   
    
    
    \subsection{Notation}
     Vectors and matrices are represented by boldface lower and upper case letters, respectively. 
  The imaginary unit is denoted by $\vec{i}$, i.e., $\vec{i}=\sqrt{-1}$.    
   The set of real, complex, and natural numbers are denoted by $\R$, $\C$, and $\N$, respectively.  
 The set of real $n$ vectors and $n{\times} m$ matrices are denoted by $\R^n$ and $\R^{n\times m}$, respectively.
 Otherwise, we use calligraphy letters to represent represent sets.
 We let $\mathcal{S}^{n-1}{=}\{\vec{x}{\in} \R^n \big|  1{=}||\vec{x}||\}$ denote the unit sphere.
 The superscript $(\cdot)\tran$ stands for transpose.
  $\diag(\vec{A}_1, {\ldots} , \vec{A}_n)$ denotes the  diagonal block matrix with $\vec{A}_1, {\ldots} , \vec{A}_n$ on the diagonal. 
 We let  $||\cdot||$ denote the $2$-norm. The distance between a point $\vec{x}\in \R^n$ and a set $\mathcal{X}\subseteq \R^n$ is given by $\dist(\vec{x},\mathcal{X})=\inf_{\vec{z}\in \R^n} || \vec{x}-\vec{z}||$. 
  For a matrix $\vec{A}\in \R^{n\times N}$, we let $\lambda_i(\vec{A})$, for $i=1,\ldots,n$ denote its eigenvalues (not in any particular order) and $\lambda_{\max}(\vec{A})=\max_{i=1,\ldots,n} |\lambda_i(\vec{A})|$. 
 We let $\lceil \vec{x}\rceil_+$ denote the projection of $\vec{x}\in \R^n$ into the positive orthant $\R_+^n$.

\section{System Model and Problem Formulation}
 \label{sec:SysModandPF}

 %
 In this section we present the power flow model, Section~\ref{subSec:2-SM}, and the voltage regulation problem,  Section~\ref{subsec:VoltReg}.
 \subsection{System Model: Linearized Power Distribution Network} \label{subSec:2-SM}

Consider a radial power distribution network with $N+1$ buses represented by the set $\mathcal{N}_0=\{0\}\cup \mathcal{N}$, where $\mathcal{N}=\{1,\ldots, N\}$.  
%
Bus $0$ is a feeder bus and the buses in $\mathcal{N}$ are branch buses. 
 Let $\mathcal{E}\subseteq \mathcal{N}_0\times \mathcal{N}_0$ denote the set of directed flow lines, so if $(i,j)\in \mathcal{E}$ then $i$ is the parent of $j$.
 For each $i$, let  $s_i=p_i+\vec{i}q_i \in \C$, $V_i\in \C$, and $v_i\in \R_+$ denote the complex power injection, complex voltage, and squared voltage magnitude, respectively,  at bus $i$.
 For each $(i,j)\in \mathcal{E}$, let $S_{ij}=P_{ij}+\vec{i}Q_{ij}\in \C$ and $z_{ij}=r_{ij}+\vec{i}x_{ij}\in \C$ denote the complex power flow and impedance in the line from bus $i$ to bus $j$.
 To model the relationship between  the variables,  we use the linearized branch flow model from~\cite{Baran1989}, which gives a good approximation in radial distribution networks.\footnote{The results also directly apply to the linearized power flow model in~\cite{Cavraro2016}. The linearized model is for the purpose of deriving the control algorithm and analyzing the theoretical performance. The designed control algorithms can be applied to the nonlinear distribution networks. In Section~\ref{Sec:Simulation}, we use the nonlinear power flow model to numerically test the algorithm performance.}
\begin{subequations} \label{eq:LinBranchFlow}
 \begin{align}
   -p_i =& P_{\sigma_i i}-\sum_{k:(i,k)\in \mathcal{E}} P_{ik},~~i\in \mathcal{N}, \\
   -q_i =& Q_{\sigma_i i}-\sum_{k:(i,k)\in \mathcal{E}} Q_{ik},~~i\in \mathcal{N}, \\
    v_j-v_i =&  -2r_{ij} P_{ij}-2 x_{ij}Q_{ij},~~(i,j)\in \mathcal{N},
 \end{align}
\end{subequations}
 where $\sigma_i$ is the parent of bus $i\in \mathcal{N}$, i.e., the unique $\sigma_i\in \mathcal{N}_0$ with $(\sigma_i,i)\in \mathcal{E}$.   
  By rearranging Equation~\eqref{eq:LinBranchFlow} we get that 
 \begin{equation}\vec{v}=\vec{A}\vec{q}^{\text{I}}+\vec{B}\vec{p}^{\text{I}}+  \vec{1}v_0,  \label{eq:PhysicalRelationship-1}\end{equation}
  where 
$\vec{v}=[v_1,\ldots,v_N]\tran$, $\vec{q}^{\text{I}}=[q_1,\ldots,q_N]\tran$, $\vec{p}^{\text{I}}=[p_1,\ldots,p_N]\tran$, 
 $$\vec{A}_{ij}=2 \hspace{-0.4cm} \sum_{(h,k)\in\mathcal{P}_i\cap \mathcal{P}_j}  \hspace{-0.4cm} x_{hk},~~\text{ and }~~\vec{B}_{ij}=2 \hspace{-0.4cm} \sum_{(h,k)\in\mathcal{P}_i\cap \mathcal{P}_j} \hspace{-0.4cm} r_{hk},$$
 where $\mathcal{P}_i\subseteq \mathcal{E}$ is the set of edges in the path from bus 0 to bus $i$.
  We use the following result in the algorithm development. 
 \begin{prop} \label{assumption:Inverse}
     $\vec{A}$ is 
      a positive definite matrix whose inverse
   has the following structure 
  \begin{equation} \label{eq:Ainverse} \displaystyle
    a_{ij}{:=}[\vec{A}^{-1}]_{ij} {=}
           \begin{cases}
              \displaystyle \frac{1}{2} \hspace{-0.05cm} \left(  x_{\sigma_i i}^{-1}{+}  \hspace{-0.25cm} \displaystyle \sum_{k:(i,k)\in\mathcal{E}}  \hspace{-0.3cm}  x_{ik}^{-1} \right) & \text{ if } i{=}j, \\
              \displaystyle -\frac{1}{2}x_{ij}^{-1}  & \text{ if } (i,j){\in} \mathcal{E} \\ 
                                                         & \text{ or } (j,i){\in} \mathcal{E}, \\
               0                      & \text{ otherwise.}
          \end{cases}
 \end{equation} 
 \end{prop}
 \begin{IEEEproof}
    It is proved in~\cite[Lemma~1]{farivar2013equilibrium} that $\vec{A}$ is positive definite.  
    Direct calculations show that $\vec{A}\vec{A}^{-1}=\vec{I}$. 
 \end{IEEEproof}

%
%

 We now introduce the Voltage Regulation Problem.

 \subsection{Voltage Regulation Problem} \label{subsec:VoltReg}

 Suppose that the real power injection $\vec{p}^{\text{I}}$ at each bus has been decided.  Write the reactive power injection as $\vec{q}^{\text{I}}=\vec{q}+\vec{q}^{\text{U}}$, where $\vec{q}$ is the adjustable reactive power that can be used for voltage regulation and $\vec{q}^{\text{U}}$ denotes other reactive power injection that cannot be changed by the voltage control devices. 
  Then the goal of the voltage regulation problem is to find feasible voltages $\vec{v}$  and adjustable reactive powers $\vec{q}$ so that the physical relationship~\eqref{eq:PhysicalRelationship-1} holds and that $\vec{v}$ and $\vec{q}$ are inside some feasible operation range $[\vec{v}^{\min},\vec{v}^{\max}]$ and $[\vec{q}^{\min},\vec{q}^{\max}]$. 
 Formally, the voltage regulation problem is to find the reactive power injection $\vec{q}$ so that,
%
\begin{subequations}   \label{eq:voltageFun-gen}
 \begin{align} \label{eq:voltageFun}
  & \vec{v}(\vec{q})= \vec{A} \vec{q}+ \vec{d},  \\
  &   \vec{v}^{\min} \leq \vec{v}(\vec{q})  \leq \vec{v}^{\max}  \label{eq:voltageCons} \\
  &   \vec{q}^{\min} \leq \vec{q} \leq \vec{q}^{\max}     \label{eq:ReactCons} 
 \end{align}  
\end{subequations} 
 where $\vec{d}:=\vec{A}\vec{q}^{\text{U}}+\vec{B}\vec{p}^{\text{I}}+\vec{1}v_0$.
  Throughout the paper we implicitly assume that the voltage regulation problem has a solution. 
 \begin{assumption}\label{assumption:feasible}
   There  exists $\vec{q}\in \R^N$  satisfying Equation~\eqref{eq:voltageFun-gen}.  
 \end{assumption}
 
 Our problem of voltage control is to design a distributed control for finding the feasible reactive power injections and voltages that satisfy equation~\eqref{eq:voltageFun-gen} that works under any operating condition denoted by $\vec{d}$.
   In particular, where each bus controls its reactive power based only  \emph{on local information and a limited number of bits communicated from neighboring buses}.   
 That is, each bus $i\in \mathcal{N}$ updates its reactive power injection according to
  \begin{equation}  \label{eq:ControlLaw}
     \vec{q}_i(t+1) = K_i(\texttt{Local\_Info}_i(t),  \texttt{Comm}_i(t)),
  \end{equation}
 where $t$ is the iteration index and 
  $K_i(\cdot)$ is the control law at bus $i$. 
  The function $K_i(\cdot)$ depends on the local information 
 $$\texttt{Local\_Info}_i(t)=(\vec{q}_i(0),\ldots,\vec{q}_i(t),\vec{v}_i(0),\ldots,\vec{v}_i(t)),$$
 and the communicated information that bus $i$ has received from its neighbors 
 $$\texttt{Comm}_i(t)=((\vec{b}_j(0))_{j\in \mathcal{N}_i},\ldots,(\vec{b}_j(t))_{j\in \mathcal{N}_i}),$$
 where $\mathcal{N}_i=\{ j \in \mathcal{N} | (i,j)\in \mathcal{E} \text{ or }  (i,j)\in \mathcal{E}\}$ and $\vec{b}_j(t)$ is the information that bus $j$ communicates to its neighbors at iteration $t$.  We study how such control laws can be achieved when $\vec{b}_i(t)$ contains a limited number of bits. 
 
 Note that the problem can generally not be solved without communication, i.e., with $\texttt{Comm}_i(t)$ empty for all $i$ and $t$, as proved by the impossibility result in~\cite{Cavraro2016}.   
 Some authors have proposed distributed solution algorithms where only  neighbors in the power networks communicate, see~\cite{Zhang_2015,vsulc2014optimal,Bolognani2013,Bolognani2015}. 
This is challenging in practice because the communication infrastructure for distribution systems is less reliable and less ubiquitous as compared to that for the bulk transmission system. To compensate for that, we study limited communication voltage control algorithms  where controllers only communicate a few bits to their neighbors. 
 These algorithms are based on combining the ideas from the voltage control algorithms in~\cite{Bolognani2013,Bolognani2015} and the limited communication gradient methods  in~\cite{Magnusson_2016}.
 However, the algorithms and results in~\cite{Magnusson_2016} are not directly applicable here because the communication structure is different and some of the assumptions in~\cite{Magnusson_2016}  do not hold for the voltage regulation problem.

\section{Voltage Control with Limited Communication} \label{Sec:algorithm}
  
 We present our limited communication algorithm for the voltage control problem in Section~\ref{Sec:Alg}. 
 We   highlight the convergence properties of the algorithm in  Section~\ref{subsec:MainConRes}.

\subsection{Algorithm} \label{Sec:Alg}


 The following algorithm is an instance of the control algorithm in Equation~\eqref{eq:ControlLaw}  where each bus communicates only 2-bits  per iteration.

\noindent \rule{\columnwidth}{2.5pt}

\noindent \textbf{VC-LB Algorithm: Voltage Control - Limited  Bandwidth}

\vspace{-0.2cm}

\noindent \rule{\columnwidth}{2.5pt}

 \begin{enumerate}[\bf a)]
 \setlength{\leftmargin}{0pt}
 \item  \textbf{Initialization:} Set $t=0$, and for each bus $i\in \mathcal{N}$ set 
\begin{align*}
 \Blambda_i(0)=&(\Blambda_i^{\min}(0),\Blambda_i^{\max}(0))=(0,0) \\
 \Bmu_i(0)=&(\Bmu_i^{\min}(0),\Bmu_i^{\max}(0))=(0,0).
\end{align*}
where the variables $\Blambda_i^{\min}(t)$, $\Blambda_i^{\max}(t)$, $\Bmu_i^{\min}(t)$, and $\Bmu_i^{\max}(t)$ are dual variables related to violating the bounds in equations~\eqref{eq:ReactCons} and~\eqref{eq:voltageCons}, see Section~\ref{SEC:CA} for details.
%

\item \textbf{Local Computation:} Each bus $i\in \mathcal{N}$ 
 computes its next reactive power injection as follows 
   \begin{align} \label{EQ:LocalQ}
     \vec{q}_i(t{+}1) {=}\Blambda_i^{\min}(t){-}\Blambda_i^{\max}(t){+}\sum_{j\in \mathcal{N}_i} a_{ij} (\Bmu_j^{\min}(t){-}\Bmu_j^{\max}(t)).
     \end{align}
  Bus $i$ can then also compute the communicated signal
 \begin{align} \label{eq:Alg1_bi}
    \vec{b}_i(t{+}1)=     \sign \left[  \begin{array}{c} \vec{q}_i(t{+}1) -\vec{q}_i^{\max} \\ \vec{q}_i^{\min}-\vec{q}_i(t{+}1) \end{array}  \right].
   \end{align}

  \item \textbf{Local Control:} Each bus $i\in \mathcal{N}$ injects the reactive power $\vec{q}_i(t{+1})$ into the power network. 
  
  \item \textbf{Local Measurement:} Each bus $i\in \mathcal{N}$ measures the voltage magnitude $\vec{v}_i(\vec{q}(t{+}1))$, given by the physical relationship~\eqref{eq:voltageFun}.

   \item \textbf{Communication:} Each bus $i\in \mathcal{N}$ communicates $\vec{b}_i(t{+}1)$ to each of its neighbours $j\in \mathcal{N}_i$ using a two bits. 

 \item \textbf{Local Computation:}      
     Each bus $i\in \mathcal{N}$ updates its dual variables
     \begin{align}
        \Blambda_i(t{+}1) {=}&  \left\lceil  \Blambda_i(t){+}\alpha \left[ 
                        \begin{array}{c} \vec{v}_i( \vec{q}(t{+}1))-\vec{v}_i^{\max}  \\ \vec{v}_i^{\min}-\vec{v}_i(\vec{q}(t{+}1))  \end{array}\right]   \right\rceil_+,    \label{eq:lamup2} \\
        \Bmu_i(t{+}1)       
                 {=}&  \left\lceil    \Bmu_i(t)   
                         {+}\beta \vec{b}_i(t) \right\rceil_+ . \label{eq:muup2}
   \end{align}
  Bus $i$ also updates a local copy of $\Bmu_j(t{+}1)$ for each neighbor $j\in\mathcal{N}_i$ using Equation~\eqref{eq:muup2}.

\item \textbf{Update Iteration Index:} $t=t+1$ and go to step \textbf{b)}. 

\end{enumerate}

\noindent \rule{\columnwidth}{2.5pt}
 
 This algorithm satisfies the structure of Equation~\eqref{eq:ControlLaw}. 
 This can be seen by noting that $\Blambda_i(t)$ is a function of $\texttt{Local\_Info}_i(t)$ via Equation~\eqref{eq:lamup2} and $\Bmu_i(t)$ is a function of 
 $\texttt{Local\_Info}_i(t)$ and $\texttt{Comm}_i(t)$ via  Equations~\eqref{eq:Alg1_bi} and~\eqref{eq:muup2}. 
 The algorithm is easy to implement since the \textbf{Local Computation} step is based on few elementary operation. 
 Moreover, for the \textbf{Communication} step each bus only needs to communicate 2 bits of information since $\vec{b}_i(t)$ can only take one of the $2^2=4$ values $(-1,-1)$, $(1,1)$, $(1,-1)$, and $(-1,1)$.  
 The parameters  $\alpha,\beta>0$ are step-sizes and are discussed further in the next section. 
 
 \begin{remark}
%
	We have considered the extreme case when only 2-bits are communicated per iterations to make the algorithm analysis more manageable. 
        It is possible that, using more bits, can improve performance and  studying such trade-offs would be an interesting extension to this work. 
	However, from both the theoretical results and numerical simulations, 2-bits communication voltage control can achieve a comparable performance as the non-quantized counterpart for our problem. 
 \end{remark}


 
 \subsection{Main Convergence Results} \label{subsec:MainConRes}
 We now study the converge of the \textbf{VC-LB Algorithm}  to a solution to the voltage regulation problem. 
 We show that in finite number of iterations the reactive power $\vec{q}(t)$ and voltage $\vec{v}(\vec{q}(t))$ satisfy Equation~\eqref{eq:voltageFun-gen}  approximately and exactly under mild additional assumptions. 
 We measure the feasibility of the reactive power $\vec{q}\in \R^N$ (and the associated voltage $\vec{v}(\vec{q})\in \R^N$) as follows 
  \begin{align} \label{eq:DistMeas} 
      \fes(\vec{q}) = \dist \big((\vec{q},\vec{v}(\vec{q})),\mathcal{Q}\times \mathcal{V}\big),
 \end{align}  
 where $\mathcal{Q}= [\vec{q}^{\min},\vec{q}^{\max}]$ and $\mathcal{V}= [\vec{v}^{\min},\vec{v}^{\max}]$.

 The following theorem (proved in Section~\ref{subsec:ProofTh1}) establishes that the \textbf{VC-LB Algorithm} can solve the problem up to any $\epsilon>0$ accuracy in finite number of iterations if the step-sizes are chosen  appropriately.  
  \begin{theorem}(Approximate Solution) \label{Theorem:DimDualGrad}
    Let $\epsilon\in(0,1]$ and 
     choose the step-sizes $\alpha,\beta>0$ such that\footnote{
     Similar convergence results 
     can be ensured for any step-size $\alpha,\beta>0$ giving in Proposition~\eqref{Lemma:Descent}
      in Section~\ref{subsec:ProofTh1}, in which case the upper bound in Equation~\eqref{eq:UB} is replaced by~\eqref{eq:UBSide}.
%
 However, we present the results here using these particular step-sizes to get an explicit form on the upper bound in~\eqref{eq:UB} in terms of available parameters. 
 }
  \begin{align} \label{eq:SSMain}
     \alpha =  \frac{1}{L}   ~~\text{ and }~~   \beta =  \frac{\epsilon}{4 N^{3/2} L}  
   \end{align} 
   where 
  $$L = \max_{i=1,\ldots,N} 2\left( \lambda_i(\vec{A}) + \frac{1}{\lambda_i(\vec{A})} \right).$$
  Then  there exists $T\in \N$ such that 
        $ \fes(\vec{q}(T)) \leq \epsilon$  
     where $T$ is upper bounded by
   \begin{align} \label{eq:UB} 
     T\leq \left\lceil  \frac{16N^3L Q \lambda_{\max}(\vec{A}) }{\epsilon^2}  \right\rceil, 
   \end{align}
   where $Q=\max \left\{(\vec{q}_i^{\min})^2,(\vec{q}_i^{\max})^2 ~:i=1,\ldots,N\right\}$.
 \end{theorem} 
 
  Theorem~\ref{Theorem:DimDualGrad} shows that  the \textbf{VC-LB Algorithm} can solve the voltage regulation problem to any precision in finite number of iterations and using finite number of communicated bits.
 In particular, the theorem shows that any  $\epsilon>0$ accuracy can be reached by communicating  $\mathcal{O}\left(1/\epsilon^2 \right)$ bits. 
 This is comparable to convergence rate without quantization. 
 This follows from the fact that the \textbf{VC-LB Algorithm} is a quantized version of a dual gradient method where the primal variable is $\vec{q}$ and the dual gradient is Lipschitz continuous, See Section~\ref{SEC:CA}. 
 Under these assumptions the feasibility of the primal iterates $\vec{q}(t)$ of dual gradient methods (with constant step-size) is $\mathcal{O}(1/\epsilon^2)$, see, for example~\cite{Beck_2014}. 
{\color{black}
 We note that to compute $L$ requires the knowledge of $\vec{A}$. However, since $\vec{A}$ is generally not changing $L$ needs to be computed once and can then be used every time when the algorithm runs.}

 In practice, we often like to run the algorithm over an extended period of time. 
 Therefore, it is desirable that once the algorithm converges to some $\epsilon>0$ accuracy then  it does not oscillate much away from that accuracy. 
 We show this numerically in Section~\ref{Sec:Simulation}. 
 However,  it is hard to prove this analytically due to technical difficulties explained in Section~\ref{SeC:Proof2}. 
 Nevertheless, the following theorem (proved in Section~\ref{SeC:Proof2}) shows that  there exist step-sizes that ensure that any $\epsilon>0$ accuracy is held for all sufficiently large $t$. 
 \begin{theorem} \label{Theorem:PD-obj-its} 
   Suppose that the voltage regulation problem is strictly feasible, i.e., there exists some $\vec{q}\in \R^N$ such that\footnote{Note that this is a stronger requirement that given by Assumption~\ref{assumption:feasible}}
    $$\vec{q}\in(\vec{q}^{\min},\vec{q}^{\max}) ~~~\text{  and }~~~ \vec{v}(\vec{q})\in(\vec{v}^{\min},\vec{v}^{\max}).$$
  Then for any $\epsilon>0$ there exist step-sizes $\alpha,\beta>0$ and $T\in \N$ such that 
\begin{align}
   \fes(\vec{q}(t)) \leq\epsilon, ~~~\text{ for all}~~ t\geq T. \label{eq:Lconv}
\end{align}
 \end{theorem}

  We can also use the \textbf{VC-LB Algorithm} to solve the voltage regulation problem exactly provided a mild additional assumption.  
    In particular, if we use the \textbf{VC-LB Algorithm} to solve the restricted voltage regulation problem  of finding reactive power $\vec{q}$ that solves
    \begin{subequations}\label{eq:ResVC}
    \begin{align} \vec{q}\in[\vec{q}^{\min}+\vec{1}\rho,\vec{q}^{\max}-\vec{1}\rho]\\  
                         \vec{v}(\vec{q})\in[\vec{v}^{\min}+\vec{1}\rho,\vec{v}^{\max}-\vec{1}\rho],
    \end{align}
  \end{subequations}
  for some $\rho>0$. 
  Then any $\epsilon=\rho$ accurate solution to the voltage regulation problem in Equation~\eqref{eq:ResVC} is an exact solution to the original problem. 
 Hence, we can apply Theorem~\ref{Theorem:DimDualGrad}   to find an exact solution in finite number of iterations where each bus communicates finite number of bits.
%
    \begin{corollary}[Exact Solution] \label{Theorem:Accurate} 
    Let $\rho >0$ be given and consider the \textbf{VC-LB Algorithm} for solving the restricted the voltage regulation problem in Equation~\eqref{eq:ResVC}.\footnote{The change in the  \textbf{VC-LB Algorithm} is that in Equation~\eqref{eq:Alg1_bi} we replace $\vec{q}_i^{\max}$ and $\vec{q}_i^{\min}$ by $\vec{q}_i^{\max}-\rho$ and $\vec{q}_i^{\min}+\rho$, respectively, and in Equation~\eqref{eq:lamup2} we replace $\vec{v}_i^{\max}$ and $\vec{v}_i^{\min}$ by $\vec{v}_i^{\max}-\rho$ and $\vec{v}_i^{\min}+\rho$, respectively.} 
    Moreover, suppose restricted the voltage regulation problem in Equation~\eqref{eq:ResVC} is feasible.  
     Choose the step-sizes $\alpha,\beta>0$ such that
  \begin{align*} 
     \alpha =  \frac{1}{L}   ~~\text{ and }~~   \beta =  \frac{\rho}{4N^{3/2} L}.
   \end{align*} 
  Then  there exists $T\in \N$ such that 
        $ \fes(\vec{q}(T)) \leq 0$  
     where $T$ is upper bounded by
   \begin{align} \label{eq:UB-2}
     T\leq \left\lceil  \frac{16N^3 LQ \lambda_{\max}(\vec{A}) }{\rho^2}  \right\rceil.
   \end{align}
 \end{corollary}
  The  parameter $\rho$  should be chosen so that the intervals in Equation~\eqref{eq:ResVC} are not empty.  
  Also note that increasing $\rho$ restricts the set of feasible voltages $\vec{v}(\vec{q}))\in[\vec{v}^{\min}+\vec{1}\rho,\vec{v}^{\max}-\vec{1}\rho]$, which must be non empty for the theorem to hold. 
 



 {\color{black}
	\subsection{Implementation without Violating Capacity Constraint}\label{subsec:projection}

	As shown in Section \ref{Sec:Alg} and \ref{subsec:MainConRes}, while the \textbf{VC-LB} algorithm can guarantee  that the reactive power $\mathbf{q}(t)$ will reach the capacity constraint $[\vec{q}^{\min},\vec{q}^{\max}]$ (approximately), $\mathbf{q}(t)$ may still violate these  limits during the transient. 
%
%
	This can be problematic since many physical control devices do not allow violating the capacity constraint for an extended amount of time. 
 One way to reduce this violation is to use a smaller  step size $\alpha$. By doing this, the voltage constraint multiplier $\Blambda(t)$ will update slower compared to the capacity constraint multiplier $\Bmu(t)$, hence putting more priority on enforcing the capacity constraint. We illustrate this in simulations, see Section~\ref{Sec:Simulation}.  
	
  We now propose an alternative implementation of the \textbf{VC-LB} algorithm that meets the reactive capacity constraint during the transients inspired by~\cite{qu2018optimal,cavraro2016value}. At time $t$, instead of implementing $\mathbf{q}_i(t)$ that may violate the capacity constraint, we implement $\mathbf{q}_i^{phy}(t) $, defined to be the projection of $\mathbf{q}_i(t)$ onto the capacity constraint set $[\mathbf{q}_i^{\min}, \mathbf{q}_i^{\max}]$. Mathematically, we implement $\mathbf{q}_i^{phy}(t) = \max( \min(\mathbf{q}_i(t),\mathbf{q}_i^{\max}),\mathbf{q}_i^{\min})$. Here the superscript in variable $\mathbf{q}_i^{phy}(t)$ refers to ``physical'', which is to emphasize that $\mathbf{q}_i^{phy}(t)$ is the amount of reactive power this is \emph{physically} implemented. Clearly, $\mathbf{q}_i^{phy}(t)$ will always satisfy the capacity constraint. We formally introduce the revised VC-LB below. 
	
	\noindent \rule{\columnwidth}{2.5pt}
	
	\noindent \textbf{VC-LB with Projection (VC-LB-P)}
	
	\vspace{-0.2cm}
	
	\noindent \rule{\columnwidth}{2.5pt}

	\begin{enumerate}[\bf a)]
		\setlength{\leftmargin}{0pt}
		\item  \textbf{Initialization:} The same as \textbf{VC-LB} Algorithm.
		%
		\item \textbf{Local Computation:} The same as \textbf{VC-LB} Algorithm.
		
		\item \textbf{Local Control:} Each bus $i\in \mathcal{N}$ injects the reactive power $\mathbf{q}_i^{phy}(t+1) = \max( \min(\mathbf{q}_i(t+1),\mathbf{q}_i^{\max}),\mathbf{q}_i^{\min})$ into the power network. 
		
		\item \textbf{Local Measurement:} Each bus $i\in \mathcal{N}$ measures the voltage magnitude $\vec{v}_i(\vec{q}^{phy}(t{+}1))$, given by the physical relationship~\eqref{eq:voltageFun}.

		\item \textbf{Communication:} The same as VC-LB Algorithm.

		\item \textbf{Local Computation:}      
		Each bus $i\in \mathcal{N}$ updates its dual variables
		\begin{align}
		\Blambda_i(t{+}1) {=}&  \left\lceil  \Blambda_i(t){+}\alpha \left[ 
		\begin{array}{c} \vec{v}_i( \vec{q}^{phy}(t{+}1))-\vec{v}_i^{\max}  \\ \vec{v}_i^{\min}-\vec{v}_i(\vec{q}^{phy}(t{+}1))  \end{array}\right]   \right\rceil_+,    \label{eq:lamup2_new} \\
		\Bmu_i(t{+}1)       
		{=}&  \left\lceil    \Bmu_i(t)   
		{+}\beta \vec{b}_i(t) \right\rceil_+ . \label{eq:muup2_new}
		\end{align}
		Bus $i$ also updates a local copy of $\Bmu_j(t{+}1)$ for each neighbor $j\in\mathcal{N}_i$ using Equation~\eqref{eq:muup2}.

		\item \textbf{Update Iteration Index:} $t=t+1$ and go to step \textbf{b)}.  
	\end{enumerate}
\noindent \rule{\columnwidth}{2.5pt}

The difference between \textbf{VC-LB-P} with the original \textbf{VC-LB} lies in step c), d), f). In step c), we implement $\mathbf{q}_i^{phy}(t+1)$ instead of $\mathbf{q}_i(t+1)$. As a result, in step d) the measured voltage becomes $\mathbf{v}_i(\mathbf{q}^{phy}(t+1))$ instead of $\mathbf{v}_i(\mathbf{q}(t+1))$, and correspondingly in step f), update equation   \eqref{eq:lamup2_new} uses $\mathbf{v}_i(\mathbf{q}^{phy}(t+1))$ instead of $\mathbf{v}_i(\mathbf{q}(t+1))$. 

We now explain the rationale behind the new implementation \textbf{VC-LB-P}. Note that in \textbf{VC-LB-P}, the update equation for Lagrangian multiplier $\Bmu_i(t)$ and variable $\vec{q}_i(t)$ is kept the same as \textbf{VC-LB}. This means that $\Bmu_i(t)$ still reflects the capacity constraint violation of $\mathbf{q}_i(t)$, so $\vec{q}_i(t)$ will still meet the capacity constraint asymptotically, and hence $\mathbf{q}_i(t) - \mathbf{q}^{phy}_i(t)\rightarrow 0$, i.e., the projection step will not substantially change the convergence results established in Section~\ref{subsec:MainConRes}. We note that similar projection steps have been done in \cite{qu2018optimal,cavraro2016value} and~\cite{qu2018optimal}  has theoretically shown that under the projection step, a modified primal-dual algorithm can still converge. Though we could not rigorously prove the convergence of VC-LB-P in this paper, we will conduct extensive simulation to verify the convergence of VC-LB-P. 
}

%
 
\section{Convergence Analysis}  \label{SEC:CA}

 The goal of this section is to prove the convergence results from the previous section. 
We give preliminary results in Section~\ref{subsec3prel} and then prove Theorems~\ref{Theorem:DimDualGrad} and~\ref{Theorem:PD-obj-its}, respectively, in Sections~\ref{subsec:ProofTh1} and~\ref{SeC:Proof2}.  
  Readers that not interested in these proofs can go straight to Section~\ref{Sec:Simulation} with no loss of information.

\subsection{Preliminaries: Duality Theory} \label{subsec3prel}

 We prove that the algorithm converges to an approximate solution to 
  the following optimization problem~\footnote{The cost can be interpenetrated  as a network loss, see \cite{Bolognani2013,Bolognani2015}.}  
 \begin{equation} \label{main_problem}
   \begin{aligned}
    & \underset{\vec{q}}{\text{minimize}}
    & & \frac{1}{2} \vec{q}\tran \vec{A} \vec{q} \\
    & \text{subject to}
    & &  \vec{v}^{\min} \leq \vec{v}(\vec{q})  \leq \vec{v}^{\max} \\
    &&&  \vec{q}^{\min} \leq \vec{q} \leq \vec{q}^{\max} .   
  \end{aligned}
\end{equation}
Problem~\eqref{main_problem} is convex because of Proposition~\ref{assumption:Inverse}. Its dual is  
\begin{equation} \label{dual_problem}
   \begin{aligned}
    & \underset{\vec{z}=(\Blambda,\Bmu)}{\text{maximize}}
    & & D(\vec{z})\\
    & \text{subject to}
    & &    \vec{z} \geq 0,
  \end{aligned}
\end{equation}
where $\Blambda =(\Blambda^{\min},\Blambda^{\max})$ and $\Bmu = (\Bmu^{\min},\Bmu^{\max})$ are the dual variables  associated to the constraints $\vec{v}^{\min} \leq \vec{v}(\vec{q})  \leq \vec{v}^{\max}$ and $ \vec{q}^{\min} \leq \vec{q} \leq \vec{q}^{\max}$, respectively, and $D:\R^{4N}\rightarrow \R$ is the dual function, see Chapter 5 in~\cite{nonlinear_bertsekas} for the details.  
 The dual gradient is
\begin{align} \label{eq:dualGrad}
  \nabla  D(\vec{z}) = \left[\begin{array}{c}  \nabla^{\Blambda}  D(\vec{z})  \\  \nabla^{\Bmu}  D(\vec{z}) \end{array}\right]
\end{align}
 where 
\begin{align*}
   \nabla^{\Blambda}  D(\vec{z}) {=} \left[\begin{array}{c} \vec{v}^{\min}{-}\vec{v}(\vec{q}(\vec{z}))  \\ \vec{v}(\vec{q}(\vec{z}) ){-}\vec{v}^{\max} \end{array}\right],~
   \nabla^{\Bmu}  D(\vec{z}) {=} \left[\begin{array}{c}  \vec{q}^{\min}{-}\vec{q}(\vec{z}) \\ \vec{q}(\vec{z}){-}\vec{q}^{\max} \end{array}\right], 
\end{align*} 
 and  
\begin{align} 
   \vec{q}(\Blambda,\Bmu) 
    =&  \Blambda^{\min}-\Blambda^{\max} + \vec{A}^{-1} \Bmu^{\min}-  \vec{A}^{-1} \Bmu^{\max}. \label{eq:LocalProblem-c} 
\end{align}
 In Lemma~\ref{Lemma:LipM} we show that $\nabla D(\cdot)$ is $L$-Lipschitz continuous, with $L$ given in the lemma. 
 Therefore, the gradient  decent method  
\begin{subequations} \label{eq:dual descent}
  \begin{align} 
     \Blambda(t{+}1) =&   \lceil \Blambda(t)+ \gamma   \nabla^{\Blambda}  D(\Blambda(t),\Bmu(t)) \rceil_+ \label{eq:dual descent_Lambda} \\
   \Bmu(t{+}1) =&  \lceil \Bmu(t) + \gamma   \nabla^{\Bmu}   D(\Blambda(t),\Bmu(t)) \rceil_+ \label{eq:dual descent_mu}
  \end{align}
\end{subequations} 
 converges to the set of optimal dual variables for appropriate step-size $\gamma$~\cite[Chapter 2]{Book_Nesterov_2004}.  
 The work in~\cite{Bolognani2013,Bolognani2015,Cavraro2016} shows how the dual gradient iterations in Equations~\eqref{eq:dual descent_Lambda} and~\eqref{eq:dual descent_mu} can be implemented in a distributed manner among the buses so that only neighbouring buses need to communicate per iteration. 
 However, these algorithms  communicate real numbers, which is challenging in practice as communication among controllers is generally constrained to low data rates.  
 To compensate for that, we have presented the \textbf{VC-LB Algorithm} in Section~\ref{Sec:algorithm}, which can be equivalently be written as  follows:
\begin{subequations} \label{eq:ALG-0}
  \begin{align}   
      \Blambda(t{+}1) =&  \lceil  \Blambda(t) +   \alpha \nabla^{\Blambda}D (\Blambda(t))  \rceil_+, \label{eq:Alg-1} \\
       \Bmu(t{+}1)       =&   \lceil \Bmu(t) +\beta \sign(\nabla^{\Bmu}D (\Bmu(t))) \rceil_+, \label{eq:Alg-2} 
 \end{align}
\end{subequations}
  where $\alpha,\beta>0$ are step-sizes and the primal variables are updated according to~\eqref{eq:LocalProblem-c}.  
 We now prove  Theorems~\ref{Theorem:DimDualGrad} and~\ref{Theorem:PD-obj-its} by considering  the \textbf{VC-LB Algorithm} on the form in Equation~\eqref{eq:ALG-0}.

 \subsection{Proof of Theorem~\ref{Theorem:DimDualGrad}} \label{subsec:ProofTh1}
  We now prove Theorem~\ref{Theorem:DimDualGrad}. 
  The main step of the proof is illustrated in the  following result (proved in Appendix~\ref{AP:Lemma:Descent}).\footnote{  Note that Proposition~\ref{Lemma:Descent} is similar to Lemma~4 in~\cite{Magnusson_2016}. 
 However, in Lemma 4 in~\cite{Magnusson_2016} the gradient is assumed to be bounded. 
 Moreover, unlike in~\cite{Magnusson_2016} the dual algorithm in this paper is an hybrid between the non-quantizes gradient step in Equation~\eqref{eq:Alg-1} and  the quantized gradient step in~\eqref{eq:Alg-2}. 
 Therefore, the results in~\cite{Magnusson_2016} do not apply here.}
  \begin{prop}  \label{Lemma:Descent} 
   Suppose $\epsilon>0$ and $\vec{z}=(\Blambda,\Bmu)\in\R_+^{4N}$ are such that $V(\vec{z}) > \epsilon$ where\footnote{Note that $V(\vec{z})=0$ is equivalent to the primal/dual variables  $(\vec{q}(\vec{z}),\vec{z})$ satisfying the KKT conditions  for optimality, see Lemma 2 in~\cite{Magnusson_2016}.}
    \begin{align} \label{eq:Lconvergence}
        V(\vec{z})= || \vec{z}-\lceil \vec{z}+ \nabla D(\vec{z}) \rceil_+||.
    \end{align}  
  Choose the step-sizes $\alpha,\beta>0$ 
\begin{subequations}  \label{eq:MainStepsize}
  \begin{align}
     \alpha <&   \frac{2}{L} \label{eq:stepsize-alpha}  \\
     \beta <& \min \left\{  \frac{\epsilon}{2N^{3/2} L},  \sqrt{\frac{\alpha(1-L \alpha/2)\epsilon^2 }{2N L}}  \right\}.\label{eq:stepsize-beta}
  \end{align} 
\end{subequations}
  Then for
      $$ \bec{z}=(\bar{\Blambda},\bar{\Bmu})= \left[   \begin{array}{c} \Blambda \\ \Bmu \end{array} +  \begin{array}{c} \alpha\nabla^{\Blambda}D (\Blambda) \\  
               \beta \sign(\nabla^{\Bmu}D (\Bmu))\end{array} \right]$$
     following holds   
   \begin{align} \label{eq:inLemmaDescent-Case0}
               D(\lceil \bar{\vec{z}}\rceil^+)    \geq D(\vec{z})+ \delta(\alpha,\beta)
   \end{align}
       where 
    \begin{align}  
    \delta(\alpha,\beta){=}& \min \left\{ \hspace{-0.1cm} \left(\alpha{-} \frac{L}{2} \alpha^2 \hspace{-0.1cm} \right)\hspace{-0.1cm} \frac{\epsilon^2}{2} {-}NL\beta^2 ,  \frac{\epsilon \beta L N }{2N^{3/2}L} {-}\beta^2 L N     \right\} \notag \\ 
     {>} &~~0 \label{eq:delta-alpha-beta}
   \end{align}
  where $L$ is an Lipschitz constant on $\nabla D(\cdot)$.
  \end{prop}

 From Lemma~\ref{Lemma:Vfun} in Appendix~\ref{APP:IMLEMMS} we have that  $\fes(\vec{q}(\vec{z}))\leq V(\vec{z})$.
 Hence, Proposition~\ref{Lemma:Descent} shows that for all $\epsilon>0$, we can choose  step-sizes $\alpha,\beta>0$ so that if $\fes(\vec{q}(\vec{z}))>\epsilon$ then the dual objective function value is improved by taking a step of the algorithm in Equations~\eqref{eq:Alg-1} and~\eqref{eq:Alg-2}.  
 Using this intuition, we have the following result. 
  \begin{prop}  \label{prop:Mainprop}
  Consider the \textbf{VC-LB Algorithm} and take some $\epsilon>0$. Choose the step-sizes $\alpha,\beta>0$ as in Equations~\eqref{eq:stepsize-alpha}  and~\eqref{eq:stepsize-beta}.  
  Then  there exists $T\in \N$ such that   $ \fes(\vec{q}(T))\leq V(\vec{z}) \leq \epsilon$    
     where $T$ is upper bounded by
   \begin{align} \label{eq:UBSide}
     T\leq \left\lceil  \frac{D^{\star}-D(\vec{z}(0))}{\delta(\alpha,\beta)}  \right\rceil.
   \end{align}
  \end{prop}
  \begin{IEEEproof}
   Suppose that $V(\vec{z}(t))>\epsilon$,  for $t=0,\ldots,T_0-1$, where $T_0:= \lceil  (D^{\star}-D(\vec{z}(0)))/\delta(\alpha,\beta) \rceil.$
  Then by  Proposition~\ref{Lemma:Descent} we have that $D(\vec{z}(t)) \geq D(\vec{z}(t{-}1))+\delta(\alpha,\beta)$ for $t=0,\ldots,T_0$ or
  \begin{align}
     0\leq& D^{\star}-D(\vec{z}(T_0))   
       \leq  D^{\star}-D(T_0-1)- \delta(\alpha,\beta)   \notag \\
                                        \leq&  D^{\star} -D(0)- T_0 \delta(\alpha,\beta) \leq0, \notag 
  \end{align}
  where the last inequality comes by that $T_0\geq (D^{\star}-D(\vec{z}(0)))/\delta(\alpha,\beta)$. 
 Hence, $\vec{z}(T_0)$ is an optimal solution to the dual problem in Equation~\eqref{dual_problem} and $\vec{q}(\vec{z}(T_0)$ is an optimal solution to the primal problem in Equation~\ref{main_problem} implying $V(\vec{z}(T_0))=0\leq \epsilon$ from Lemma~\ref{Lemma:Vfun}. 
  \end{IEEEproof}
  Theorem~\ref{Theorem:DimDualGrad} follows from Proposition~\ref{prop:Mainprop}. 
  Direct calculation shows that the step-sizes in Equation~\eqref{eq:SSMain} satisfy the condition given in Equations~\eqref{eq:stepsize-alpha}  and~\eqref{eq:stepsize-beta}. Therefore, we just have to show for $\alpha$ and $\beta$ in Equation~\eqref{eq:SSMain} following holds 
  \begin{equation} \label{eq:LemmaNUB}
      \frac{D^{\star}-D(\vec{z}(0))}{\delta(\alpha,\beta)} \leq  \frac{16N^3L  Q  \lambda_{\max}(\vec{A}) }{\epsilon^2}
   \end{equation}
   to conclude the proof. 
    Equation~\eqref{eq:LemmaNUB} follows directly from the following two identity (proved below):
   \begin{align*}
       \delta(\alpha,\beta) =& \frac{3\epsilon^2}{16 N^2 L}, ~~~\text{ and} ~~~
       D^{\star}-D(\vec{z}(0))\leq  \lambda_{\max}(\vec{A}) N Q, 
   \end{align*}
   for all $\hec{q}\in [\vec{q}^{\min},\vec{q}^{\max}]$ such that $\vec{v}(\hec{q})\in [\vec{v}^{\min},\vec{v}^{\max}]$.  
  The equality is obtained  by  plugging $\alpha$ and $\beta$ from Equation~\eqref{eq:SSMain} into Equation~\eqref{eq:delta-alpha-beta}.  
   To prove the inequality we note that $D^{\star}=(\vec{q}^{\star})\tran \vec{A}\vec{q}^{\star}$ because of the strong duality and 
$D(\vec{z}(0))=\min_{\vec{q}\in\R^N} \vec{q}\tran \vec{A}\vec{q}= 0$, since $\vec{z}(0)=\vec{0}$. 
 Yielding  
 $$ D^{\star}{-}D(\vec{z}(0))\leq (\vec{q}^{\star})\tran \vec{A}\vec{q}^{\star}\leq \lambda_{\max}(\vec{A}) ||\vec{q}^{\star}||^2\leq \lambda_{\max}(\vec{A}) N Q.$$

 \subsection{Proof of Theorem~\ref{Theorem:PD-obj-its}}\label{SeC:Proof2} 

  We now prove Theorem~\ref{Theorem:PD-obj-its}.  
  The result is based on the following result proved in Appendix.
 \begin{prop}  \label{Lemma:boundedZ}
   Suppose that the voltage regulation problem is strictly feasible, i.e., there exists some $\vec{q}\in \R^N$ such that 
    $$\vec{q}\in(\vec{q}^{\min},\vec{q}^{\max}) ~~~\text{  and }~~~ \vec{v}(\vec{q})\in(\vec{v}^{\min},\vec{v}^{\max}).$$
   Define the set  $\mathcal{Z}(\epsilon) = \left\{ \vec{z}\in \R_+^{4N} | V(\vec{z})\leq \epsilon \right\}$ and the function  $F(\epsilon)=\inf\{ D(\vec{z}) | \vec{z}\in\mathcal{Z}(\epsilon)\}$.
   Then following holds:
\begin{enumerate}[a)]
   \item There exists $\kappa>0$ such that the set $\mathcal{Z}(\epsilon)$ is bounded and $F(\epsilon)>-\infty$ for all $\epsilon\leq \kappa$. 

  \item We have  $\lim_{\epsilon\rightarrow 0^+}F(\epsilon)=f^{\star}$.

  \item           For any $\epsilon>0$, step-sizes $\alpha,\beta>0$ chosen as in Equations~\eqref{eq:MainStepsize} and $T\in \N$ such that $V(\vec{z}(T))\leq \epsilon$ we have 
   \begin{align} \label{eq:Dlemma}
     D(\vec{z}(t))\geq F(\epsilon)- \frac{L}{2}\beta^2, \text{ for all } t\geq T.
   \end{align} 
\end{enumerate}
 \end{prop}

 From Proposition~\ref{prop:Mainprop} and Proposition~\ref{Lemma:boundedZ} we have that for any $\eta>0$  there exit step-sizes $\alpha,\beta>0$ and $T\in \N$ such that $D^{\star}- D(\vec{z}(t)) \leq \eta$ holds for all $t\geq T$.
 From Lemma~\ref{Lemma:LevelSets} in Appendix~\ref{APP:IMLEMMS}, 
 there exists $\eta>0$ such that if $D^{\star}-D(\vec{z})\leq \eta$ then $\dist(\vec{z},\mathcal{Z}^{\star})\leq \epsilon/\phi$, where $\phi=||\vec{G}||$ and
 $$\vec{G}= \left( \begin{array}{cccc} \vec{I}   & -\vec{I}    & \vec{A}^{-1} & - \vec{A}^{-1} \\
                                               \vec{A}  & - \vec{A}  & \vec{I}   & -\vec{I} 
                  \end{array} \right).$$ 
 Therefore, there exists $\alpha,\beta>0$ and $T\in \N$ such that  $\dist(\vec{z}(t),\mathcal{Z}^{\star})\leq\epsilon/\phi$ for all $t\geq T$. 
 We also have for  $\vec{z}\in \R^{4N}$,  by setting $\vec{z}^{\star}= \min_{\bec{z}\in\mathcal{Z}^{\star}} ||\vec{z}-\bec{z}||$, that 
\begin{align*}
   \texttt{fes}(\vec{q}(\vec{z})) \leq&  \big|\big| (\vec{q}(\vec{z}),\vec{v}(\vec{q}(\vec{z})))-(\vec{q}(\vec{z}^{\star}),\vec{v}(\vec{q}(\vec{z}^{\star}))\big|\big|,  \\
   =&||\vec{G}(\vec{z}-\vec{z}^{\star})||\leq \phi ||\vec{z}-\vec{z}^{\star}||= \phi \dist(\vec{z},\mathcal{Z}^{\star}),
\end{align*}
 where we have used that  $(\vec{q}(\vec{z}),\vec{v}(\vec{q}(\vec{z})))=\vec{G}\vec{z}+(\vec{0},\vec{d})$. 
%
 Therefore,  $ \texttt{fes}(\vec{q}(\vec{z}(t))) =\dist(\vec{z}(t),\mathcal{Z}^{\star})\leq\epsilon$ for all $t\geq T$.

\section{Numerical Results using Non Linear Power Flow} \label{Sec:Simulation}

\begin{figure}
    	\centering
    	\includegraphics[width=0.48\textwidth]{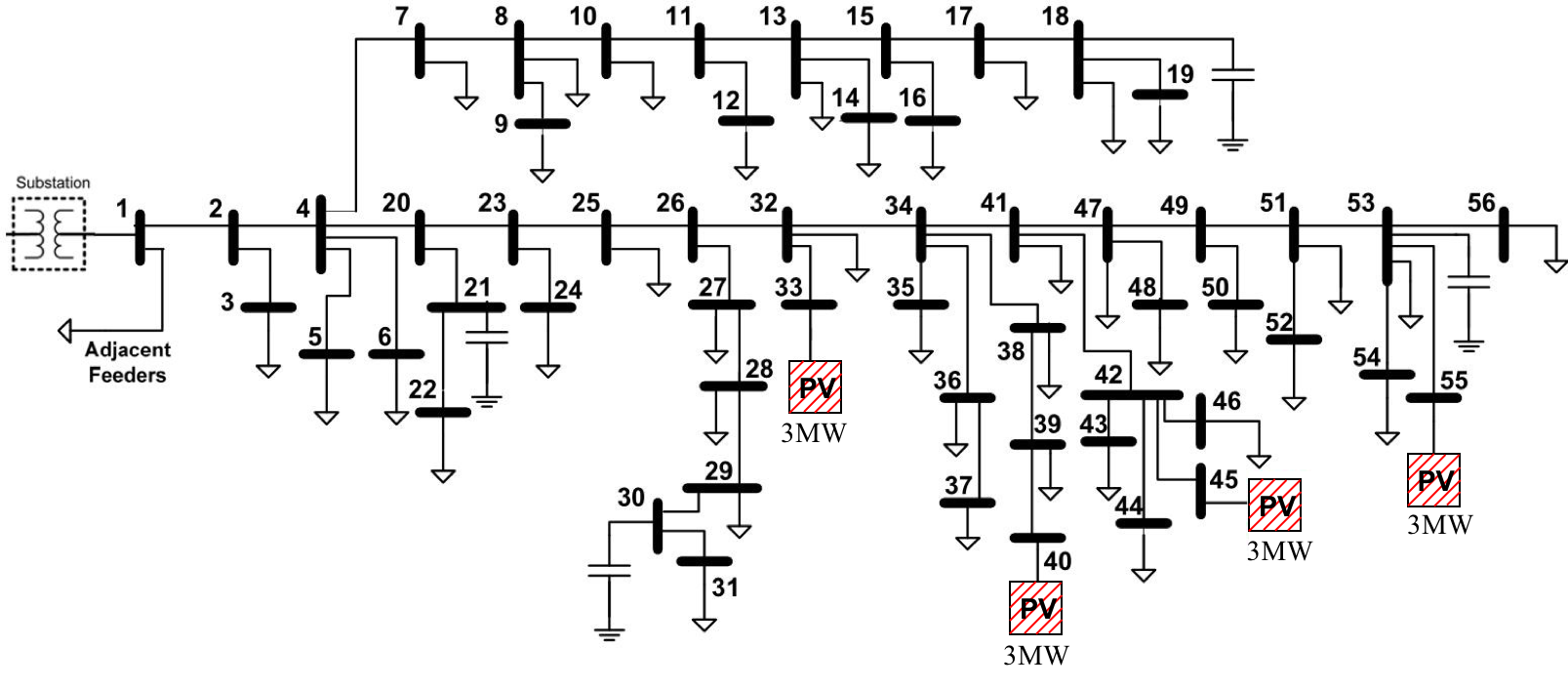}
    	\caption{The 56 bus network used in our simulation. }\label{fig:56bus}
    \end{figure}

   We test our algorithm on the 56-bus radial distribution network in~\cite{farivar2012optimal}, as shown in Figure \ref{fig:56bus}. Bus 1 is the feeder bus, and there are PVs installed at bus 33, 40, 45, 55. All quantities, when units are not explicitly given, are in the per unit (p.u) system. The nominal value of voltage is $12$kV. Throughout the simulation we set $v^{\min}=11.4$kV, $v^{\max}=12.6$kV ($\pm 5\%$ of the nominal value), and for each bus $i\in \mathcal{N}$, we set $\vec{q}_i^{U}=0$, $\vec{q}_i^{\min}=-0.5$MW, and $\vec{q}_i^{\max}=0.5$MW. The power flow is calculated using MATPOWER \cite{zimmerman2011matpower},\footnote{All simulations are run on MATLAB R2016b on Macbook Pro 2015 Model with 2.7 GHz Intel Core i5.} which uses the full nonlinear power flow instead of the linearized power flow in our analysis. We demonstrate the proposed algorithm in static and dynamic voltage control scenarios. 
    
    \subsection{Static Voltage Control}
    
\begin{figure}
    \centering

        \includegraphics[width =\columnwidth]{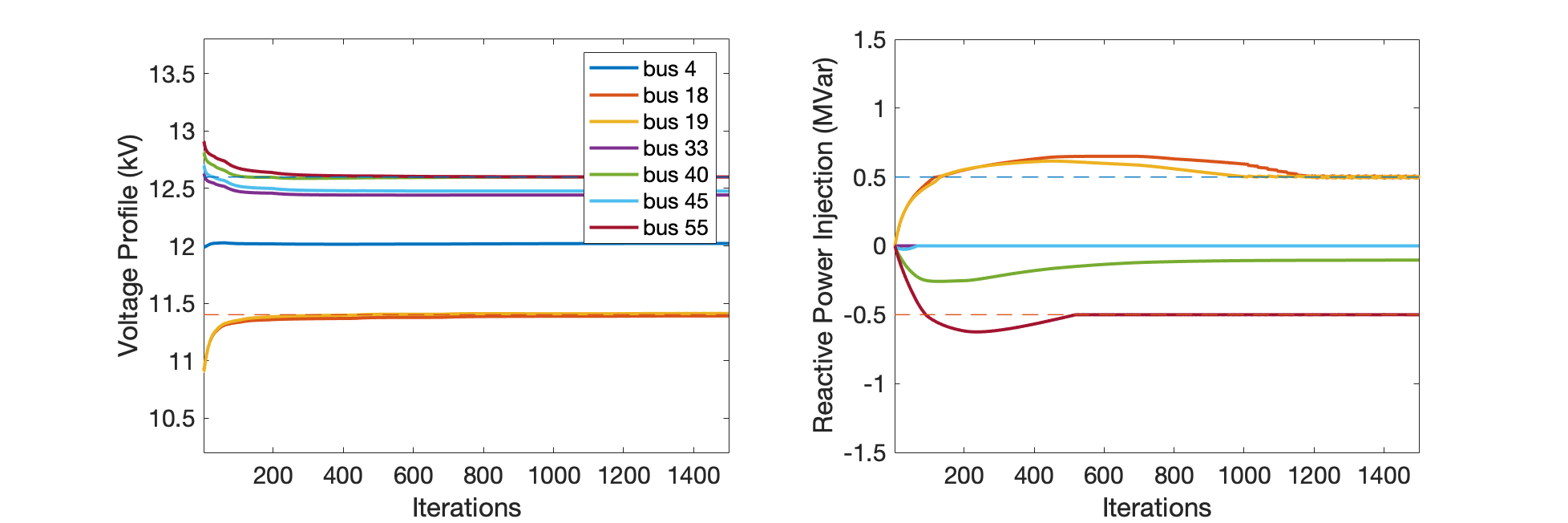}

    \caption{Convergence of the $\vec{v}(t)$ and $\vec{q}(t)$ in the \textbf{VC-LB Algorithm}. Parameters are $\alpha = 0.2,\beta = 10^{-5}, \rho = 0$. Notice that during each iteration, each node uses 2-bits to communicate.}\label{fig:Dist}
    \label{fig:STATCon}
\end{figure}~

\begin{figure}
    \centering
        \includegraphics[width=\columnwidth]{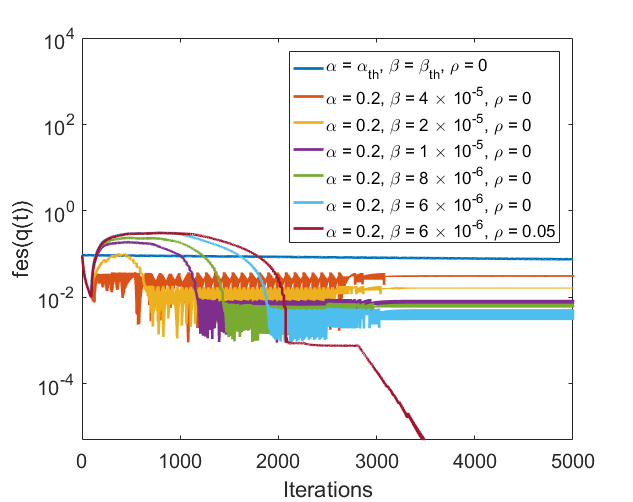}
    \caption{$\texttt{fes}(\vec{q}(t))$ in Equation~\eqref{eq:DistMeas} (in per unit system) under different values of $\alpha$, $\beta$ and $\rho$. $\alpha_{\text{th}}$ and $\beta_{\text{th}}$ are the theoretic value of step sizes given in Equation~\eqref{eq:SSMain} with $\epsilon
    	 = 1$ and their values are $\alpha_{\text{th}} = 1.4093\times 10^{-4}$, $\beta_{\text{th}} = 6.4061\times 10^{-7}$ under per unit system. The $\alpha=0.2,\beta = 1\times 10^{-5}$, $\rho = 0$ case has the final values of $\texttt{fes}(\vec{q}(t))$ dropping below $1\times 10^{-2}$, and the corresponding trajectory of $\vec{v}(t)$ and $\vec{q}(t)$ are in Fig. ~\ref{fig:STATCon}. }
 \label{fig:FES}
\end{figure}
We first test the case in which the real power injection $\mathbf{p}_i$ is fixed and is from the load data in \cite{farivar2012optimal}, except we have scaled up the load from bus 7 to bus 19 to make the problem more challenging. The parameters of the algorithm is set as $\alpha=0.2$, $\beta=10^{-5}$ and $\rho = 0$. 
    Fig.~\ref{fig:STATCon} illustrates the voltage profile and the reactive power injection trajectory under \textbf{VC-LB}.
   The results show that the voltage converge to the feasible range within roughly 400 iterations, or 800 bits of communication.  
  Similarly, the reactive power converges to the feasible range within  roughly 1200 iterations, or 2400 bits of communication.

  Fig.~\ref{fig:FES}  demonstrates the convergence of \textbf{VC-LB} under different step-sizes $\alpha,\beta$ and parameter $\rho$ (recall $\rho$ is introduced to find an exact solution, see Corollary \ref{Theorem:Accurate}).
  It plots the feasibility measure $\text{fes}(\vec{q}(t))$ in Equation~\eqref{eq:DistMeas} as a function of iteration counts.
 First of all, Fig.~\ref{fig:FES} shows that when $\rho = 0$, the algorithm converges to an approximately feasible point, within certain accuracy. Further, there is a trade-off between how high the accuracy is and how many communication bits are needed to achieve that accuracy. Secondly, for the $\rho > 0$ case, Fig.~\ref{fig:FES} demonstrates our algorithm is indeed driving the system towards an exact solution, consistent with Corollary \ref{Theorem:Accurate}.
 

{\color{black}We comment that the step-sizes in Equation~\eqref{eq:SSMain} are conservative. For the parameters of this network, we can calculate $L\approx 7\times10^3$, so the step-sizes in Equation~\eqref{eq:SSMain} are very small ($\alpha \approx 10^{-4}$, $\beta \approx 6\times 10^{-7}$ using $\epsilon=1$), leading to very slow convergence, as illustrated in the blue curve in Fig.~\ref{fig:FES}. On the contrary, the results in  Fig.~\ref{fig:STATCon} and Fig.~\ref{fig:FES} have shown that in practice fast convergence can be obtained for step-sizes that are much larger than \eqref{eq:SSMain}. 
 This is because to obtain the theoretical guarantees we must account for every potential worst-case behavior and also make some relaxations to make the mathematical derivations tractable. This is a typical trade-off between theoretical convergence guarantees and the convergence in practice, as discussed for example in~\cite{nonlinear_bertsekas}.}

\begin{figure}
	\centering
	
	\includegraphics[width =\columnwidth]{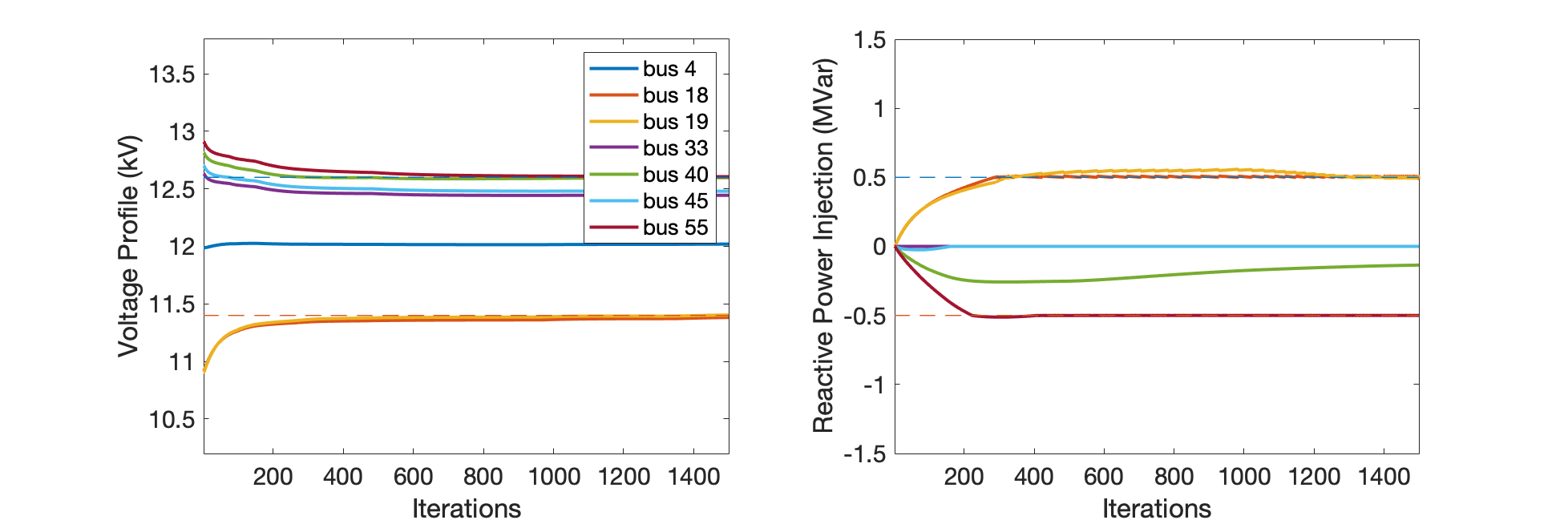}
	
	\caption{Convergence of the $\vec{v}(t)$ and $\vec{q}(t)$ in the \textbf{VC-LB} Algorithm. Parameters are $\alpha = 0.08,\beta = 10^{-5}, \rho = 0$. }
	\label{fig:smallalpha}
\end{figure}

\begin{figure}
	\centering
	
	\includegraphics[width =\columnwidth]{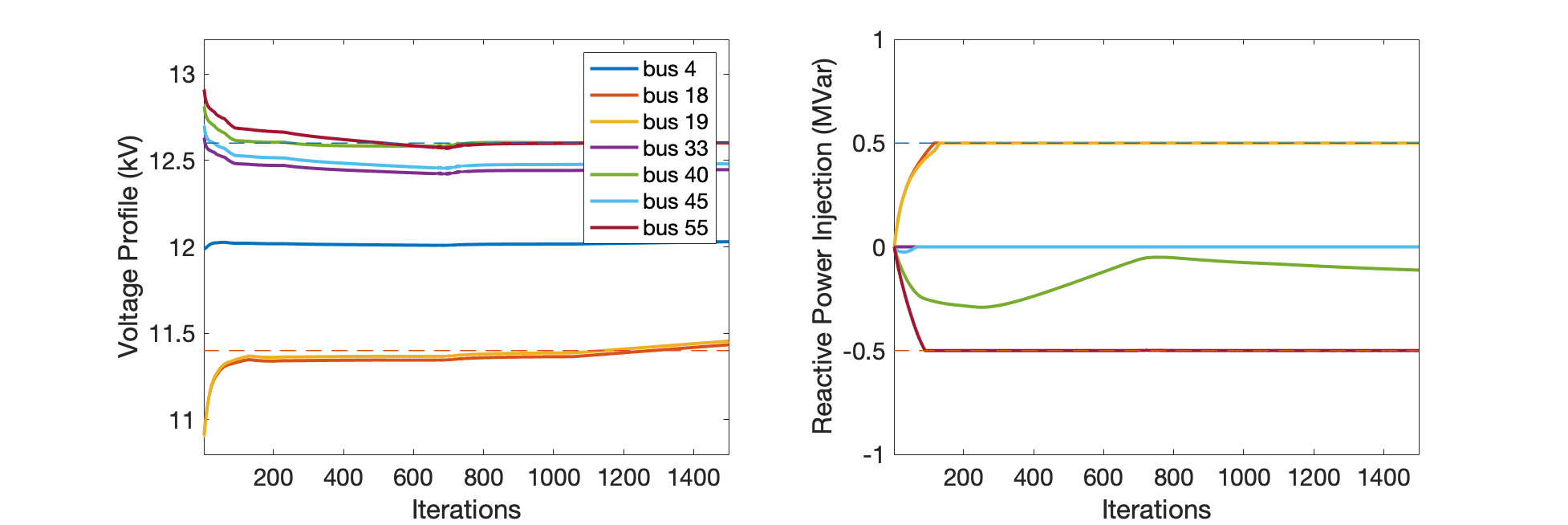}
	
	\caption{Convergence of the $\vec{v}(t)$ and $\vec{q}^{phy}(t)$ in the \textbf{VC-LB-P} Algorithm. Parameters are $\alpha = 0.2,\beta = 10^{-5}, \rho = 0$. Note that the reactive power injection never exceeds the upper and lower limit. }
	\label{fig:projection}
\end{figure}

{\color{black} We note that in Fig. \ref{fig:STATCon}, though the reactive power injection meets the capacity constraint asymptotically, it violates the constraint during the transient. As discussed in Section~\ref{subsec:projection}, we now test the two methods to deal with this. Firstly, we test the \textbf{VC-LB} algorithm with a small $\alpha$, which as discussed in Section~\ref{subsec:projection}, will put more priority on enforcing capacity constraint. We use the same simulation setting and step sizes as Fig.~\ref{fig:STATCon}, except that $\alpha$ is reduced to $0.08$ ($\beta=10^{-5}$, $\rho=0$ are kept the same). Results are shown in Fig.~\ref{fig:smallalpha}. Compared to Fig.~\ref{fig:STATCon}, the capacity constraint violation is now only minor in Fig.~\ref{fig:smallalpha}. Interestingly, we note that compared to Fig.~\ref{fig:STATCon}, the voltage profile in Fig.~\ref{fig:smallalpha} converges slower, which makes sense since in Fig.~\ref{fig:smallalpha} we have put more priority on enforcing capacity constraint as opposed to the voltage constraint. Secondly, we test the variant of \textbf{VC-LB} algorithm, the \textbf{VC-LB-P} algorithm in Section~\ref{subsec:projection}. Simulation settings are the same as
%
%
 Fig.~\ref{fig:STATCon} and we use the same step sizes $\alpha=0.2,\beta=10^{-5},\rho=0$. The voltage profile and the reactive power injection profile are presented in Fig.~\ref{fig:projection}. Fig.~\ref{fig:projection} confirms that the \textbf{VC-LB-P} algorithm still converges and in the meanwhile the reactive power injection does not violate the capacity constraint at any time.   }
 
 ~

    \subsection{Dynamic Voltage Control}

\begin{figure}[t!]
    \centering
        \includegraphics[width=\columnwidth]{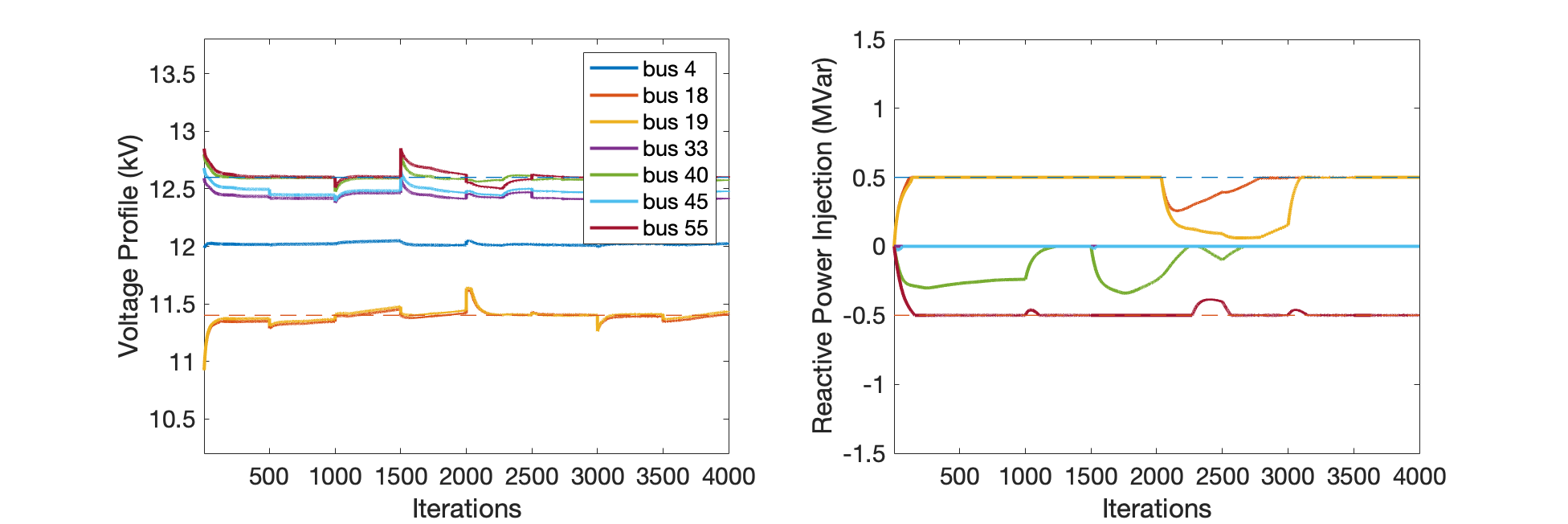}


    \caption{Dynamic voltage control: convergence of $\vec{v}(t)$ and $\vec{q}^{phy}(t)$ in the \textbf{VC-LB-P Algorithm}. Parameters are $\alpha=0.2,\beta=10^{-5},\rho=0$. }\label{fig:dynamic}
    \label{fig:DYNCon}
\end{figure}

{\color{black}
 In practice voltage regulation algorithms must respond quickly to fluctuating electric behaviors of the consumers of the network. To this end,
  we test the  \textbf{VC-LB-P} Algorithm in a dynamic environment where the real power injections change over the course of the algorithm.  
 Specifically, we let the real power injections fluctuate (randomly) to a new set of values once in a while, and within each interval between two fluctuations, the algorithm \textbf{VC-LB-P Algorithm}  is allowed to run $500$ iterations (communicate $1000$ bits) to respond to the fluctuation.
The value of the fluctuating real power injection is determined through multiplying the static value of the real power injection (the value used in the previous subsection) by a random scalar drawn uniformly from $[0.75,1.25]$. The parameters of \textbf{VC-LB-P} Algorithm is set as $\alpha=0.2$, $\beta=10^{-5}$, and $\rho = 0$, the same as that of Fig.~\ref{fig:projection}. 
 
 The voltage profile and the reactive power injection of a selected set of buses over $8$ intervals ($4000$ iterations) are shown in Fig.~\ref{fig:DYNCon}. 
 The figure shows that in all cases the \textbf{VC-LB-P} Algorithm regulates the voltage profile within the 500-iteration (1000-bit) limit and in the meanwhile not violating the reactive power capacity constraint at any time. We further note that the time it takes to conduct a $1$-bit communication, when using protocols for extremely low latency communications~\cite{Durisi_2016} (at the cost of limited data-rate), can be made extremely small. 
This suggests that the \textbf{VC-LB-P} Algorithm in this paper can be used to provide fast voltage regulation in future smart grids. 
}

\section{Conclusion}

This paper studied distributed voltage control algorithms where only a few bits of communication between neighboring buses are needed. 
 The convergence of these algorithms  was studied and their practical applicability illustrated in simulations.  
 Future work is to 1) study  the trade-offs between the data rate and control performance, 2) characterize the online performance under dynamic operating scenarios,  and 3) implement the algorithm on a real-world testbed.


\bibliographystyle{IEEEbib}
\bibliography{refs}



\appendices

\section{Important Lemmas} \label{APP:IMLEMMS}
\begin{lemma}  \label{Lemma:LipM}
  We have $\nabla D(\vec{z})=\vec{M}\vec{z}+\vec{r}$, where 
  \begin{align} \label{eq:matrixM}
 \vec{M}= \left( \begin{array}{cccc} -\vec{A} & \vec{A} &  -\vec{I} & \vec{I} \\ 
                                     \vec{A} & -\vec{A} &  \vec{I} & -\vec{I} \\
                                    -\vec{I} & \vec{I} &  -\vec{A}^{-1} & \vec{A}^{-1}\\
                                    \vec{I} & -\vec{I} &  \vec{A}^{-1} & -\vec{A}^{-1}
           \end{array} \right).
 \end{align}
 The eigenvalues of $\vec{M}$ are
 \begin{align} \label{eq:LambdaM} 
       \lambda_i(\vec{M})=-2 \left( \lambda_i(\vec{A}) + \frac{1}{\lambda_i(\vec{A})}\right),~~\text{ for } i=1,\ldots,n
 \end{align}
 and $\lambda_i(\vec{M})=0$ for $i=n+1,\ldots,4n$. 
%
\end{lemma}
\begin{IEEEproof}
 The fact that $\nabla D(\vec{z})=\vec{M}\vec{z}+\vec{r}$ follows from Equations~\eqref{eq:dualGrad} and~\eqref{eq:LocalProblem-c}.  
 To find the eigenvalues of $\vec{M}$, let $\vec{w}_i\in \R^n$ be an eigenvector of the matrix $\vec{A}$ associated to the eiginvalue  $\lambda_i(\vec{A})$. 
 Then direct calculations show that $\vec{M}\vec{v}_i = \lambda_i(\vec{M})\vec{v}$, with $\lambda_i(\vec{M})$ defined in Equation~\eqref{eq:LambdaM}, where 
 $\vec{v}_i=(\lambda_i(\vec{A}) \vec{w}_i,-\lambda_i(\vec{A}) \vec{w}_i,  \vec{w}_i, -\vec{w}_i)$.

 To prove that $0$ is a an eigenvalue of $\vec{M}$ with multiplicity 3n, i.e., $\lambda_i(\vec{M})=0$ for $i=n+1,\ldots,4n$, we note that $\vec{M}$ is a rank $n$ matrix and hence has a $3n$ dimensional null space.  
 This can be seen by noting that the rows 2, 3, and 4 of the block matrix $\vec{M}$, see Equation~\eqref{eq:matrixM}, can be obtained from elementary row operations (Gauss elimination) on the first row of $\vec{M}$.  In particular, $\vec{R}_2=-\vec{R}_1$, $\vec{R}_3=\vec{A} \vec{R}_1$, $\vec{R}_4=- \vec{A}\vec{R}_1$, where $\vec{R}_i$ is the i-th row of the block representation of $\vec{M}$ given in Equation~\eqref{eq:matrixM}.  
\end{IEEEproof}
 \begin{lemma}  \label{Lemma:Vfun}
 For $V(\cdot)$ defined in Equation~\eqref{eq:Lconvergence} we have:
 \begin{enumerate} [a)]
   \item A feasible dual variable $\vec{z}\in \R_+^{4N}$ is an optimal solution to the Dual Problem~\eqref{dual_problem} if and only if $V(\vec{z})=0$. 

   \item For all $\vec{z}\in \R_+^{4N}$ 
        $\fes(\vec{q}(\vec{z}))\leq V(\vec{z})$, 
   where $\vec{q}(\vec{z})$ is defined in Equation~\eqref{eq:LocalProblem-c}.


 \end{enumerate}

 \end{lemma} 
 \begin{IEEEproof}
  a) Follows from~\cite[Lemma~3]{Magnusson_2016}.

   b) We have that
 $$ \fes(\vec{q}(\vec{z})){=} 
       \left|\left| \left\lceil \left[ \begin{array}{c} 
       \vec{v}^{\min}-\vec{v}(\vec{q}(\vec{z}))  \\ 
       \vec{v}(\vec{q}(\vec{z}) )-\vec{v}^{\max} \\ 
       \vec{q}^{\min}-\vec{q}(\vec{z})           \\ 
       \vec{q}(\vec{z})-\vec{q}^{\max}          
       \end{array}\right] \right\rceil_+ \right| \right| 
  {=} || \lceil \nabla D(\vec{z}) \rceil_+||, 
$$
 where the later equality comes from Equation~\eqref{eq:dualGrad}. 
  Using that 
 $$ \lceil \nabla_i D(\vec{z}) \rceil_+ \leq |\vec{z}_i-\lceil \vec{z}_i + \nabla_i D(\vec{z}) \rceil_+ |$$
  for $i=1,\ldots,4N$, we get that 
 \begin{multline*}
        || \lceil \nabla D(\vec{z}) \rceil_+||^2= \sum_{i=1}^{4N}  \lceil \nabla_i D(\vec{z}) \rceil_+^2 \\
           \leq 
                   \sum_{i=1}^{4N} |\vec{z}_i-\lceil \vec{z}_i + \nabla_i D(\vec{z}) \rceil_+ |^2=V(\vec{z})^2,
 \end{multline*}
 which proves the result.
 \end{IEEEproof}
 
   \begin{lemma}\label{Lemma:projections}
 For all $\beta \in [0,1]$, $z\in \R$ and $x,\alpha_1,\alpha_2 \in \R_+$  with $\alpha_1 \leq |x-\lceil x+ z \rceil_+|$ following holds 
 \begin{align}
  \beta |x-\lceil x+z\rceil^+ | \leq& |x- \lceil x+ \beta z\rceil^+| \label{eq:LemPro-1}\\
    \alpha_1 =& |x-\lceil x+\alpha_1 \sign(z) \rceil_+ | \label{eq:LemPro-2} \\
   0 \leq& z(\lceil x+\alpha_2 z\rceil^+-x ). \label{eq:LemPro-3} 
 \end{align}
\end{lemma}
\begin{IEEEproof}
 The proof follows similar steps as the proof of \cite[Lemma~9]{Magnusson_2016}.
 In particular, the inequalities~\eqref{eq:LemPro-1} and~\eqref{eq:LemPro-2} are easily checked by using  that for $x\in \R_+$ and $z\in \R$
 we have $|x-\lceil x+ z\rceil^+| =|z|$  if $x+z\geq 0$ and $|x-\lceil x+ z\rceil^+| =x$  if $x+z\leq 0$.
%
 The Inequality~\eqref{eq:LemPro-3} comes by rearrange the inequality 
 $\sign(z)\lceil x+\alpha_1 z \rceil_+\leq  \sign(z)x$ and multiply $|z|$ on both sides.  
\end{IEEEproof}

   \begin{lemma}\label{Lemma:optbounded} 
   Suppose that the voltage regulation problem is strictly feasible, i.e., there exists some $\vec{q}\in \R^N$ such that $\vec{q}\in(\vec{q}^{\min},\vec{q}^{\max})$ and $ \vec{v}(\vec{q})\in(\vec{v}^{\min},\vec{v}^{\max})$. 
  Then the set $\mathcal{Z}^{\star}$ of optimizers of the dual problem given in Equation~\ref{dual_problem} is bounded. 
\end{lemma}
\begin{IEEEproof} 
  Follows directly from Lemma~1 in~\cite{nedic2009approximate}. 
\end{IEEEproof}

\begin{lemma}   \label{Lemma:LevelSets}
  For $\epsilon>0$, there exists $\eta>0$ such that if $D^{\star}-D(\vec{z})\leq \eta$ then $\dist(\vec{z},\mathcal{Z}^{\star})\leq \epsilon$.  
\end{lemma}
\begin{IEEEproof}
 We prove this by contradiction. 
   Suppose that no such $\eta>0$ exists.  
   Then we can generate a sequence $\vec{z}_k\in \R_+$ such that $\lim_{k\rightarrow \infty}D(\vec{z}_k)=D^{\star}$ and $\dist(\vec{z}_k,\mathcal{Z}^{\star})>\epsilon$ for all $k\in \N$.
   The sequence $\vec{z}_k$ is bounded because $\mathcal{Z}^{\star}$ is bounded, see Lemma~\ref{Lemma:optbounded}. 
  Therefore, the level sets 
 $\mathcal{L}(\kappa)= \{\vec{z}\in \R_+^{4N} | D(\vec{z})\geq D^{\star}-\kappa  \}$ 
are also bounded, see~\cite[Proposition B.9]{nonlinear_bertsekas}. 
   As a result, $\vec{z}_k$ has a convergent subsequence $\vec{z}_{k_r}$ with  $\lim_{r\rightarrow \infty } \vec{z}_{k_r} =\bar{\vec{z}}\in \R_+$. 
   Then $D(\bar{\vec{z}})=\lim_{r\rightarrow \infty} D(  \vec{z}_{k_r} )=D^{\star}$, so $\bar{\vec{z}}\in \mathcal{Z}^{\star}$.
   This contradicts the fact that $\dist(\vec{z}_k,\mathcal{Z}^{\star})>\epsilon$ for all $k\in \N$, since $\dist(\vec{z}_{k_r},\mathcal{Z}^{\star}) \leq ||\vec{z}_{k_r}-\bar{\vec{z}}||\leq \epsilon$ for sufficiently large $r$, which yields the result. 
\end{IEEEproof}

\section{Proof of Proposition~\ref{Lemma:Descent}} \label{AP:Lemma:Descent}
  The dual function $D(\cdot)$ is concave and with $L$-Lipschitz continues gradient so we have~\cite[Equation (2.1.6)]{Book_Nesterov_2004}
   \begin{align}
    D(\lceil \bar{\vec{z}}\rceil^+) {\geq} & D(\vec{z}) {+} \langle \nabla D(\vec{z}), \lceil \bar{\vec{z}}\rceil^+{-}\vec{z}    \rangle  
      {-} \frac{L}{2} ||\lceil \bar{\vec{z}} \rceil^+{-} \vec{z}||^2, \notag \\
            {=} & D(\vec{z}) +\langle \nabla^{\Blambda} D(\vec{z}),  \lceil \bar{\Blambda} \rceil_+ {-}\Blambda \rangle - \frac{L}{2} ||\lceil \bar{\Blambda} \rceil_+-\Blambda||^2 \notag \\
              &~~~~~~ +\langle \nabla^{\Bmu} D(\vec{z}),  \lceil \bar{\Bmu} \rceil_+ {-}\Bmu \rangle - \frac{L}{2} ||  \lceil \bar{\Bmu}\rceil_+ - \Bmu ||^2 \notag \\
             \geq& D(\vec{z}) + \alpha\left(1-\frac{L}{2} \alpha \right) || \Blambda -\lceil \Blambda + \nabla^{\Blambda} D(\vec{z})\rceil_+ ||^2 \notag \\ 
            & {+}\langle \nabla^{\Bmu} D(\vec{z}), \lceil \bar{\Bmu} \rceil_+{-} \Bmu \rangle {-} \frac{L}{2} ||  \lceil \bar{\Bmu} \rceil_+ {-} \Bmu||^2, \label{eq:inLemmaDescent}
\end{align}     
  where the final inequality is proved below. 
 From Lemma~\ref{Lemma:Vfun} we have $V(\vec{z})\geq \texttt{fes}(\vec{q}(\vec{z}))>\epsilon$, yielding
 \begin{align*} \epsilon^2 \leq& V(\vec{z})^2  
                                    {=} || \Blambda{-}\lceil  \Blambda{+} \nabla^{\Blambda} D(\Blambda)\rceil_+||^2{+}|| \Bmu{-}\lceil  \Bmu{+} \nabla^{\Bmu} D(\Bmu)\rceil_+||^2.
 \end{align*}
 Consider the two cases separately when (a) 
 $$|| \Bmu{-}\lceil  \Bmu{+} \nabla^{\Bmu} D(\Bmu)\rceil_+||^2\geq \frac{\epsilon^2}{2}$$
 and (b) \vspace{-0.3cm} 
  \begin{align*}
     || \Bmu{-}\lceil  \Bmu{+} \nabla^{\Bmu} D(\Bmu)\rceil_+||^2 <& \frac{\epsilon^2}{2} ~~\text{ and } \\
       || \Blambda-\lceil  \Blambda{+} \nabla^{\Blambda} D(\Blambda)\rceil_+||^2\geq& \frac{\epsilon^2}{2}.
  \end{align*}
 In case (a) we have  (proved below)
  \begin{align} \label{eq:inLemmaDescent-CaseA}
       D(\lceil \bar{\vec{z}}\rceil^+)     \geq D(\vec{z})  + \left( \frac{\epsilon }{2N^{3/2} L} - \beta \right) \beta L N,
  \end{align}
  and case (b)  we have (proved below)
   \begin{align}
    D(\lceil \bar{\vec{z}}\rceil^+)         \geq & D(\vec{z}) + \alpha\left(1-\frac{L}{2} \alpha \right)\frac{\epsilon^2}{2}  -  NL\beta^2  \label{eq:inLemmaDescent-CaseB}.
\end{align} 
  Then Equation~\eqref{eq:inLemmaDescent-Case0} is obtained by Equations~\eqref{eq:inLemmaDescent-CaseA} and~\eqref{eq:inLemmaDescent-CaseB}. 
      
    \underline{\textbf{Proof of Equation~\eqref{eq:inLemmaDescent}:}} 
    We have 
    \begin{align}
      \langle \nabla^{\Blambda} D(\vec{z}),  \lceil \bar{\Blambda} \rceil_+ {-}\Blambda \rangle {=}&
              \sum_{i=1}^{2N} |\nabla_i^{\Blambda} D(\vec{z})| | \lceil \Blambda_i {+} \alpha \nabla_i^{\Blambda} D(\vec{z})\rceil_+ {-}\Blambda_i|    \notag    \\
              {\geq}&  \sum_{i=1}^{2N} \frac{1}{\alpha}    (\lceil \Blambda_i {+} \alpha \nabla_i^{\Blambda} D(\vec{z})\rceil_+ {-} \Blambda_i)^2 \label{eq:DescLemma-eq2}
    \end{align}
    where the equality comes by using that each term of the sum is positive from Equation~\eqref{eq:LemPro-3} and the inequality comes by using the non-expansiveness of the projection to get \vspace{-0.2cm} 
  $$ |\nabla_i^{\Blambda} D(\vec{z})| \geq \frac{1}{\alpha} \left|  \lceil \Blambda_i +  \alpha \nabla_i^{\Blambda} D(\vec{z})\rceil_+ - \Blambda_i\right|. \vspace{-0.1cm} $$
  We also have \vspace{-0.3cm} 
  \begin{align}
      || \lceil \bar{\Blambda} \rceil_+ -\Blambda||^2 = \sum_{i=1}^{2N} (\lceil \Blambda_i {+} \alpha \nabla_i^{\Blambda} D(\vec{z})\rceil_+-\Blambda_i)^2.  \label{eq:DescLemma-eq3}
  \end{align}
  By combining~\eqref{eq:DescLemma-eq2} and~\eqref{eq:DescLemma-eq3} we get
    \begin{align*}
       \langle \nabla^{\Blambda} D(\vec{z}), & \lceil  \bar{\Blambda} \rceil_+ -\Blambda \rangle - \frac{L}{2} ||  \lceil \bar{\Blambda} \rceil_+ - \Blambda||^2 \\
          \geq& \left(\frac{1}{\alpha} - \frac{L}{2}\right)\sum_{i=1}^{2N}( \lceil \Blambda_i {+} \alpha \nabla_i^{\Blambda}  D(\vec{z})\rceil_+ - \Blambda_i)^2 \\
          \geq& \alpha^2 \left(\frac{1}{\alpha} - \frac{L}{2}\right)\sum_{i=1}^{2N}( \lceil \Blambda_i {+}  \nabla_i^{\Blambda}  D(\vec{z})\rceil_+ - \Blambda_i)^2 \\
             \geq& \alpha \left(1-\frac{L}{2}\alpha \right)||\Blambda {-} \lceil \Blambda {+}  \nabla^{\Blambda} D(\vec{z})\rceil_+||^2,   
    \end{align*}
    where the second inequality comes by using~\eqref{eq:LemPro-1} and the fact that $\alpha<2/L<1$ (see Lemma~\ref{Lemma:LipM}).

    \underline{\textbf{Proof of Equation~\eqref{eq:inLemmaDescent-CaseA}:}} Suppose that
 $$||\Bmu - \lceil \Bmu + \nabla^{\Bmu} D (\Bmu)   \rceil_+||^2\geq  \epsilon^2/2.$$ 
  Then we have
    \begin{align}
       \langle  \nabla^{\Bmu} &D(\vec{z}),  \lceil \bar{\Bmu} \rceil_+ - \Bmu \rangle \notag  \\
      {=}&  \sum_{i=1}^{2N} |\nabla_i^{\Bmu} D(\vec{z})| | \lceil \Bmu_i {+} \beta \sign( \nabla_i^{\Bmu}  D(\vec{z}))\rceil_+-\Bmu_i| \notag \\
         \geq &   |\Bmu_{j}{-}\lceil \Bmu_{j} {+}  \nabla_{j}^{\Bmu} D(\vec{z})\rceil_+| |\Bmu_{j}{-} \lceil \Bmu_{j} {+} \beta\sign( \nabla_{j}^{\Bmu} D(\vec{z}))\rceil_+ | \notag \\ 
  \geq &   \frac{\epsilon}{2 \sqrt{N}} \beta \label{eq:mainLemmaINEQ2-a}
    \end{align}
  where $j={\text{argmax}}_{i=1,\cdots,2N}  |\Bmu_i{-} \lceil \Bmu_i {+}  \nabla_i^{\Bmu} D(\vec{z})\rceil_+|$,  
  the equality comes by using that each term of the sum is positive from Equation~\eqref{eq:LemPro-3}, 
  the first inequality from the non-expansiveness of the projection~\cite[Proposition B.11]{nonlinear_bertsekas}, 
  and the final inequality comes from the fact that $||\Bmu - \lceil \Bmu + \nabla^{\Bmu} D (\vec{z})   \rceil_+||_{\infty}\geq  \epsilon/(\sqrt{2}\sqrt{2N})$, since $||\Bmu - \lceil \Bmu + \nabla^{\Bmu} D (\vec{z})   \rceil_+||\geq  \epsilon/\sqrt{2}$, and Equation~\eqref{eq:LemPro-2} together with the fact that 
 $\beta\leq \epsilon /(2\sqrt{N}) \leq |\Bmu_{j}{-} \lceil \Bmu_{j} {+}  \nabla_{j}^{\Bmu} D(\vec{z})\rceil_+|.$
 We also have
 \begin{align} 
     || \Bmu {-}\lceil \bar{\Bmu} \rceil_+||^2 {=}& \sum_{i=1}^{2N} |\Bmu_{i}{-} \lceil \Bmu_{i} {+} \beta\sign( \nabla_{i}^{\Bmu} D(\vec{z}))\rceil_+|^2 \notag \\ 
    \leq & 2N \beta^2. \label{eq:mainLemmaINEQ2-b}
 \end{align}
 Combining Equations~\eqref{eq:mainLemmaINEQ2-a} and~\eqref{eq:mainLemmaINEQ2-b} and rearranging yields
    \begin{equation*} 
       \langle \nabla^{\Bmu} D(\vec{z}),  \lceil \bar{\Bmu} \rceil_+ {-} \Bmu \rangle - \frac{L}{2} || \lceil \bar{\Bmu} \rceil_+ {-} \Bmu||^2  
          \geq \left(\frac{\epsilon }{2 N^{3/2} L} {-}  \beta \right) \beta L N,
    \end{equation*} 
 which proves Equation~\eqref{eq:inLemmaDescent-CaseA}.

    \underline{\textbf{Proof of Equation~\eqref{eq:inLemmaDescent-CaseB}:}} 
 The result follows from that \vspace{-0.2cm} 
    \begin{multline} \label{eq:mainLemmaIneq3} 
       \langle \nabla^{\Bmu} D(\vec{z}),  \lceil \bar{\Bmu} \rceil_+ {-}\Bmu\rangle {-} \frac{L}{2} || \lceil \bar{\Bmu} \rceil_+- \Bmu||^2 
          \geq - N L \beta^2 ,
    \end{multline}  
   where inequality comes by that    \vspace{-0.1cm} 
  $$|| \Bmu -\lceil \bar{\Bmu} \rceil_+||^2=\sum_{i=1}^{2N} (\Bmu -\lceil \bar{\Bmu} \rceil_+)^2 \leq 2N\beta^2.$$   
  and  $\langle \nabla^{\Bmu} D(\vec{z}), \lceil \bar{\Bmu} \rceil_+ {-} \Bmu \rangle\geq 0$, from Equation~\eqref{eq:LemPro-3}.

  \section{Proof of Lemma~\ref{Lemma:boundedZ}} \label{App:Proof}
   We prove the result by contradiction.
  Suppose that such an $\epsilon>0$ does not exist. 
  Then we can construct a sequence $\vec{z}^{k}  \in \R_+^N$ such that $\lim_{k\rightarrow \infty}||\vec{z}^k||=\infty$ and  $\lim_{k\rightarrow \infty}V(\vec{z}^k)=0$.   
  By considering the set
\begin{align} \label{eq:finallemma-Iset}
  \mathcal{I}=\{i=1,\cdots,4N \big| \lim_{k\rightarrow \infty} \vec{z}_i^k=0\},
 \end{align}    
 we can further restrict the sequence $\vec{z}^k$ so that for $i\notin \mathcal{I}$ it holds that $\vec{z}_i^k\geq W$, for some $W>0$ and all $k\in \N$. 

 To obtain the contradiction we consider the sequence $\bec{z}^k$ given by $\bec{z}_i^k=0$ if $i\in \mathcal{I}$ and $\bec{z}_i^k=\vec{z}_i^k$. 
 In the sequel we show the contradicting results that
 \begin{equation} \label{eq:inBoundLemma-con1} 
      \lim_{k\rightarrow \infty} V(\bec{z}^k)=0
\end{equation}
 and that there exists $\delta>0$ and $K\in\N$ such that
 \begin{equation} \label{eq:inBoundLemma-con2}
    V(\bec{z}^k)\geq \delta~~ \text{ for all } ~~k\geq K.
 \end{equation}
  Equations~\eqref{eq:inBoundLemma-con1} and~\eqref{eq:inBoundLemma-con2} clearly contradict each other so there can not exist such sequence $\vec{z}^k$, which yields the result. 
  We now prove Equations~\eqref{eq:inBoundLemma-con1} and~\eqref{eq:inBoundLemma-con2}.
  
  \vspace{0.2cm}

\noindent \underline{\textbf{Proof of Equation~\eqref{eq:inBoundLemma-con1}:}} 
 We have 
  $$ V(\bec{z}^k)^2 =\sum_{i=1}^{4N} (\bec{z}_i^k - \lceil \bec{z}_i^k+\nabla_i D(\bec{z}^k) \rceil_+)^2, $$   
 so it suffices to show that
 $$ \lim_{k\rightarrow \infty} |\bec{z}_i^k - \lceil \bec{z}_i^k+\nabla_i D(\bec{z}^k) \rceil_+|=0, ~\text{ for }~ i=1,\ldots,4N. $$ 
 Consider first the case when $i\notin \mathcal{I}$. 
 Then 
 \begin{align*}
    |\nabla_i D(\bec{z}^k)|{-}|\nabla_i D(\vec{z}^k)| 
          \leq& |\nabla_i D(\bec{z}^k)-\nabla_i D(\vec{z}^k)|\\
          \leq& L||\bec{z}^k-\vec{z}^k||= L \sqrt{\sum_{j\in \mathcal{I}} (\vec{z}_j^k)^2}
 \end{align*} 
  where the first inequality comes by the triangle inequality, the second inequality from that $\nabla D$ is $L$-Lipschitz continuous, and the equality by the definition of $\bec{z}$.  
  By rearranging, we have 
 \begin{equation} \label{eq:inLastLemmaDbz}
    |\nabla_i D(\bec{z}^k)| \leq |\nabla_i D(\vec{z}^k)| + L \sqrt{\sum_{j\in \mathcal{I}} (\vec{z}_j^k)^2},
 \end{equation}
 where the right hand side converges to zero since $$\lim_{k\rightarrow \infty} |\nabla_i D(\vec{z}^k)|=0$$ since $i\notin \mathcal{I}$, see Claim 2 below, and $\lim_{k\rightarrow \infty}\vec{z}_j^k=0$ for $j\in\mathcal{I}$. 
 Therefore,  $\lim_{k\rightarrow \infty}|\nabla_i D(\bec{z}^k)|=0$ and by using that $|\bec{z}_i^k - \lceil \bec{z}_i^k+\nabla_i D(\bec{z}^k) \rceil_+| \leq |\nabla_i D(\bec{z}^k)|$, see~\cite[Proposition B.11.(c)]{nonlinear_bertsekas}, we have
 $$ \lim_{k\rightarrow \infty} |\bec{z}_i^k - \lceil \bec{z}_i^k+\nabla_i D(\bec{z}^k) \rceil_+|=0. $$
 
 Consider next the case when $i\in \mathcal{I}$.
 We obtain Equation~\eqref{eq:inLastLemmaDbz} in this case as well, following the same steps as before. 
 By the limit of both sides of Equation~\eqref{eq:inLastLemmaDbz} we get
 \begin{multline} \label{eqInMlemma-limsup}
   \limsup_{k\rightarrow \infty} \nabla_i D(\bec{z}^k) \leq \limsup_{k\rightarrow \infty} \nabla_i D(\vec{z}^k)  \\ + \lim_{k\rightarrow \infty}L \sqrt{\sum_{j\in \mathcal{I}} (\vec{z}_j^k)^2}\leq 0,
 \end{multline}
 where the first inequality comes from~\cite[Proposition A.4.(d)]{nonlinear_bertsekas} and the second inequality comes from
 that $\limsup_{k\rightarrow \infty}\nabla_i D(\vec{z}^k)\leq 0$, see Claim 1 below, and that $\lim_{k \rightarrow \infty }\vec{z}_j^k=0$ for all $j\in \mathcal{I}$. 
 This together with that $\bec{z}_i^k=0$ yields
\begin{align*}
  \lim_{k\rightarrow \infty} |\bec{z}_i^k {-} \lceil \bec{z}_i^k{+}\nabla_i D(\bec{z}^k) \rceil_+| 
   {=}\lim_{k\rightarrow \infty} \min\{0, \nabla_i D(\bec{z})\}=0.
\end{align*}

   \vspace{0.2cm}

\noindent \underline{\textbf{Proof of Equation~\eqref{eq:inBoundLemma-con2}:}}
 Consider the sequence 
 $$ \vec{w}^k= \frac{\vec{z}^{\star}-\bec{z}^k}{||\vec{z}^{\star}-\bec{z}^k||},$$
  where $\vec{z}^{\star}$ is some element of $\mathcal{Z}^{\star}$. 
  We start by showing that there exists $K_0\in \N$ and $\kappa$ such that 
 $\langle \nabla D(\bec{z}^k),\vec{w}^k\rangle\geq \kappa$ for all $k\geq K_0$. 
 The set $\mathcal{Z}^{\star}$ is bounded by Slater's condition and Lemma~1 in~\cite{nedic2009approximate}. 
Therefore, there exists $R>0$ such that $\mathcal{Z}^{\star}\subseteq \{ \vec{z}\in \R_+^{4N}| R>||\vec{z}-\vec{z}^{\star}||\}$.  
 Then there exists $\phi>0$ such that
 \begin{equation}\label{eqInFinalLemmaPHI} 
     \phi=\min_{\vec{z}\in\mathcal{S}} \langle \nabla D(\vec{z}),\vec{z}^{\star}-\vec{z}  \rangle,
 \end{equation}
 where $\mathcal{S}=\{\vec{z}\in\R_+ | R=||\vec{z}^{\star}-\vec{z}||\}$, since $\mathcal{S}$ is compact set, the intersection $\mathcal{S}\cap \mathcal{Z}^{\star}$ is empty, and  $\langle \nabla D(\vec{z}),\vec{z}^{\star}-\vec{z}  \rangle>0$ for all $\vec{z}\in \R_+\setminus \mathcal{Z}^{\star}$, see Claim 3 below.
 Moreover, since $\lim_{k\rightarrow \infty} ||\bec{z}^k||=\infty$, there exists $K_0\in \N$ such that $ ||\bec{z}^k-\vec{z}^{\star}||> R$ for all $k\geq K_0$. 
 Therefore, we get
 \begin{align}
   \langle \nabla D(\bec{z}^k),\vec{w}^k \rangle \geq&  \langle \nabla D (\vec{z}^{\star}-R\vec{w}^k),\vec{w}^k \rangle \notag \\ 
                   =&  \frac{1}{R} \langle \nabla D (\vec{z}^{\star}-R\vec{w}^k),\vec{z}^{\star}-(\vec{z}^{\star}-R\vec{w}^k) \rangle \notag \\
                   \geq& \frac{\phi}{R}=:\kappa,
 \end{align}
 for all $k\geq K_0$, where the first inequality comes by using that $\bec{z}^k=\vec{z}^{\star}-||\vec{z}^{\star}-\bec{z}^k||\vec{w}^k$ together with the fact that $\nabla D$ is monotone decreasing, since $D$ is concave, to obtain
 \begin{align*}
    -\langle \nabla D(\bec{z}^k)-\nabla D(\vec{z}^{\star}-R \vec{w}^k),(R-||\vec{z}^{\star}-\bec{z}^k||) \vec{w}^k \rangle \geq 0,
 \end{align*}
 or by rearranging and using that $(R-||\vec{z}^{\star}-\bec{z}^k||)<0$ 
 \begin{align*}
    \langle \nabla D(\bec{z}^k), \vec{w}^k\rangle \geq \langle \nabla D(\vec{z}^{\star}-R \vec{w}^k), \vec{w}^k \rangle,
 \end{align*}
  and the final inequality comes from the fact that $\vec{z}^{\star}-R\vec{w}^k\in\mathcal{S}$ and Equation~\eqref{eqInFinalLemmaPHI}.

 From above we have that 
 \begin{align}
   \kappa \leq& \langle \nabla D(\bec{z}^k),\vec{w}^k \rangle \notag \\
     =& \langle \nabla_{\mathcal{I}} D(\bec{z}^k),\vec{w}_{\mathcal{I}}^k \rangle + \langle \nabla_{\mathcal{I}^C} D(\bec{z}^k),\vec{w}_{\mathcal{I}^C}^k \rangle, \label{eqInFinalLemma-kappabound}
 \end{align} 
 for $k\geq K_0$, where $\mathcal{I}^C=\{1,\ldots,4N \}\setminus\mathcal{I}$. 
 We also have  
\begin{align*}
 \limsup_{k\rightarrow \infty} \langle \nabla_{ \mathcal{I}} D(\bec{z}^k),\vec{w}_{\mathcal{I}}^k \rangle =&  \limsup_{k\rightarrow \infty} \sum_{i\in \mathcal{I}} \nabla_i D(\bec{z}^k)\vec{w}_i^{k}  \\
  \leq & \sum_{i\in \mathcal{I}}  \limsup_{k\rightarrow \infty} \nabla_i D(\bec{z}^k)\vec{w}_i^k  \\ 
     \leq& 0 
\end{align*}
 where the first inequality comes from~\cite[Proposition A.4.(d)]{nonlinear_bertsekas} and the second inequality comes from 
 that $\limsup_{k\rightarrow \infty}\nabla_i D(\bec{z}^k)\leq 0$, see Equation~\eqref{eqInMlemma-limsup}, and the fact that $\vec{w}_i^{k}\geq 0$ for all $k\in \N$, since $\bec{z}_i^k=0$, and $\vec{w}_i^{k}\leq 1$.
 Therefore, we can choose $K\in \N$, with $K\geq K_0$, so that $\langle \nabla_{ \mathcal{I}} D(\bec{z}^k),\vec{w}_{\mathcal{I}}^k \rangle\leq \kappa/2$ for all $k\geq K$. 
  Then from Equation~\eqref{eqInFinalLemma-kappabound} we have 
 $$ \frac{\kappa}{2}\leq \langle \nabla_{\mathcal{I}^C} D(\bec{z}^k),\vec{w}_{\mathcal{I}^C}^k \rangle~~ \text{ for all }~~ k\geq K, $$
 and by the Cauchy-Schwarz inequality we have 
 $$ \frac{\kappa}{2}\leq ||\nabla_{\mathcal{I}^C} D(\bec{z}^k)||~||\vec{w}_{\mathcal{I}^C} ||\leq ||\nabla_{ \mathcal{I}^{C}} D(\bec{z}^k)||,$$
 for all $k\geq K$.  
 By using the equivalence between the 2- and $\infty$-norms we have 
 \begin{align} \label{eqInLastLemmaGradDB} \frac{\kappa}{2\sqrt{4N}}\leq ||\nabla_{ \mathcal{I}^{C}} D(\bec{z}^k)||_\infty = |\nabla_{r_k}  D(\bec{z}^k)|,
 \end{align}
  for all $k\geq K$, where 
  $${r_k}=\underset{j\in\mathcal{I}^C}{\text{argmax}} |\nabla_j  D(\bec{z}^k)|.$$
 Hence, we get
 \begin{align*}
    V(\bec{z}^k) =&  || \bec{z}^k - \lceil \bec{z}^k+\nabla D(\bec{z}^k) \rceil_+|| \\
                \geq& |\bec{z}_{r_k}^k - \lceil \bec{z}_{r_k}^k+\nabla_{r_k} D(\bec{z}^k) \rceil_+| \\
              \geq& \delta:=\min\left\{W,\frac{\kappa}{2\sqrt{4N}}\right\},
 \end{align*}  
 for all $k\geq K$, 
 where the final inequality comes from that ${r_k}\in \mathcal{I}^C$ and Equation~\eqref{eqInLastLemmaGradDB} so $\bec{z}_{r_k}^k\geq W$ and $|\nabla_{r_k} D (\bec{z})|\geq \kappa/(2\sqrt{4N})$ for all $k\geq K$.

\vspace{0.2cm} 
 
\noindent \underline{\textbf{Claim 1:} If $i\in \mathcal{I}$ then $\limsup_{k\rightarrow \infty} \nabla_i D(\vec{z}^k)\leq 0$.} 
    Suppose the contrary, that $\limsup_{k\rightarrow \infty} \nabla_i D(\vec{z}^k)> 0$ for $i\in \mathcal{I}$. 
   Then there exists a scalar $\kappa>0$ and a subsequence $\vec{z}^{k_r}$ such that $\nabla_i D(\vec{z}^{k_r})>\kappa$ for all $r\in \N$.
   Since $\vec{z}_i^{k_r}$ and $\nabla_i D(\vec{z}^{k_r})$ are nonnegative we have 
  $$\lceil \vec{z}_i^{k_r}+ \nabla_i D(\vec{z}^{k_r}) \rceil_+= \vec{z}^{k_r}+\nabla_i D(\vec{z}^{k_r}).$$
 Therefore, it holds that
  \begin{align*}
          \kappa <& \nabla_i D(\vec{z}^{k_r}) 
                      =    |\vec{z}_i^{k_r}- \lceil \vec{z}_i^{k_r}+ \nabla_i D(\vec{z}^{k_r}) \rceil_+| \\
                     \leq& ||\vec{z}^{k_r}- \lceil \vec{z}^{k_r}+ \nabla D(\vec{z}^{k_r}) \rceil_+||_{\infty} \\
                     \leq& ||\vec{z}^{k_r}- \lceil \vec{z}^{k_r}+ \nabla D(\vec{z}^{k_r}) \rceil_+|| =V(\vec{z}^{k_r}),              
  \end{align*}
  which  contradicts the fact that $\lim_{k\rightarrow \infty} V(\vec{z}^{k})=0$.

\vspace{0.2cm} 
 
\noindent  \underline{\textbf{Claim 2:} If $i\notin \mathcal{I}$ then $\lim_{k\rightarrow \infty} \nabla_i D(\vec{z}^k)= 0$.} 
 
\noindent Suppose that $i\notin \mathcal{I}$.  Then  
 $$|\vec{z}_i^k-\lceil \vec{z}_i^k+\nabla_i D(\vec{z}^k) \rceil_+| \geq \min \{ W, |\nabla_i D(\vec{z}^k)| \},$$ 
  for all $k\in \N$ since $\vec{z}_i^k\geq W$.
 We also have that 
 $$ \lim_{k\rightarrow \infty}|\vec{z}_i^k-\lceil \vec{z}^k+\nabla D(\vec{z}_i^k) \rceil_+| =0$$
 since $\lim_{k\rightarrow \infty}V(\vec{z}^k)=0$.  
 Therefore, we have
  $$\lim_{k\rightarrow \infty} \min \{ W, |\nabla_i D(\vec{z}^k)| \}=0.$$
 Since $W>0$ it must hold that  $\lim_{k\rightarrow \infty}  |\nabla_i D(\vec{z}^k)|=0$.

\noindent  \underline{\textbf{Claim 3:} $\langle \nabla D(\vec{z}),\vec{z}^{\star}-\vec{z}\rangle >0$ for all $\vec{z}\in \R_+^{4N}\setminus \mathcal{Z}^{\star}$.} 
           %

\noindent  Note that by the KKT optimality conditions we have that $\vec{z}^{\star}\in \mathcal{Z}^{\star}$ if and only if $\vec{z}_i^{\star}\geq0$,  $\nabla_i D(\vec{z}^{\star})\leq 0$, and $\nabla_i D(\vec{z}^{\star}) \vec{z}_i^{\star}=0$ for all $i=1,\ldots,4N$. 
 Using this fact, we  now prove the result for all $\vec{z}\in \R_+^{4N}\setminus \mathcal{Z}^{\star}$ by consider the following two cases (a)  there exists $r\in\{1\ldots,4N\}$ such that $\vec{z}_r>0$ and $\nabla_r D(\vec{z}^{\star})<0$ and (b)    $\vec{z}_i=0$ for all $i$ such that $\nabla_i D(\vec{z}^{\star})<0$. 

Consider first Case (a). 
  Then we have  
  \begin{align*}  
    \langle \nabla D(\vec{z}), \vec{z}^{\star}{-}\vec{z}\rangle\geq&  \langle \nabla D(\vec{z}^{\star}), \vec{z}^{\star}{-}\vec{z}\rangle \\
         \geq&  \sum_{i=1}^N {-}\nabla_i D(\vec{z}^{\star}) \vec{z}_i\geq  -\nabla_r D(\vec{z}^{\star}) \vec{z}_r   >0,
  \end{align*}
 where the first inequality comes from that $-\nabla D$ is monotone, since $D$ is concave, the second  inequality from that $\nabla_i D(\vec{z}^{\star}) \vec{z}_i^{\star}=0$,
 the third inequality from that $\nabla_i D(\vec{z}^{\star})\leq 0$ and $\vec{z}_i\geq 0$ for all $i=1,\ldots,4N$, and
    the final inequality from the fact that $\nabla_r D(\vec{z}^{\star})<0$ and that $\vec{z}_r>0$. 

 Consider next Case (b). 
 Then   $\vec{z}\in \R_+^{4N}\setminus \mathcal{Z}^{\star}$ and  $\vec{z}_i=0$ for all $i$ such that $\nabla_i D(\vec{z}^{\star})<0$. 
 Therefore, we have $\nabla D(\vec{z})\neq D(\vec{z}^{\star})$ because if $\nabla D(\vec{z})= D(\vec{z}^{\star})$ then the optimality condition  $\vec{z}_i\geq0$,  $\nabla_i D(\vec{z})\leq 0$, and $\nabla_i D(\vec{z}^{\star}) \vec{z}_i=0$, hold for all $i=1,\ldots,4N$, so $\vec{z}\in \mathcal{Z}^{\star}$. 
  Then by using that the function $-D$ is convex with $L$-Lipschitz continuous gradient we have~\cite[eq.~(2.1.8)]{Book_Nesterov_2004} 
   \begin{multline*}
     \langle \nabla D(\vec{z}), \vec{z}^{\star}{-}\vec{z}\rangle\geq  \langle \nabla D(\vec{z}^{\star}), \vec{z}^{\star}{-}\vec{z}\rangle \\
    {+} \frac{1}{L} ||\nabla D(\vec{z}^{\star}){-}\nabla D(\vec{z}) ||^2{>}0,
   \end{multline*}
  where the final inequality comes by using that $\langle \nabla D(\vec{z}^{\star}), \vec{z}^{\star}-\vec{z}\rangle\geq 0$ from~\cite[Proposition 2.1.2]{nonlinear_bertsekas} and that $||\nabla D(\vec{z}^{\star})-\nabla D(\vec{z}) ||>0$.  

\begin{IEEEbiography}[{\includegraphics[width=1in,height=1.25in,clip,keepaspectratio]{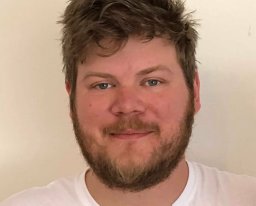}}]{Sindri Magn\'usson}
received the B.Sc. degree in Mathematics from University of Iceland, Reykjav\'ik Iceland, in 2011, the Masters degree in Mathematics from KTH Royal Institute of Technology, Stockholm Sweden, in 2013, and the PhD in Electrical Engineering from the same institution, in 2017. He is a Postdoctoral Fellow  in the School of Engineering and Applied Sciences in Harvard University. 
 His research interests include distributed optimization, both theory and applications.
\end{IEEEbiography}

 \begin{IEEEbiography}[{\includegraphics[width=1in,height=1.25in,clip,keepaspectratio]{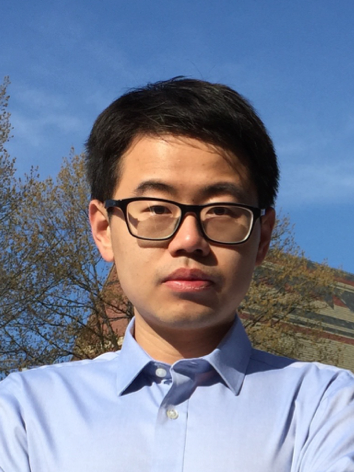}}]{Guannan Qu} 
 received his B.S. degree in Electrical Engineering from Tsinghua University in Beijing, China in 2014. Since 2014 he has been a graduate student in the School of Engineering and Applied Sciences at Harvard University. His research interest lies in network control and optimization.
 \end{IEEEbiography}

 \begin{IEEEbiography}[{\includegraphics[width=1in,height=1.25in,clip,keepaspectratio]{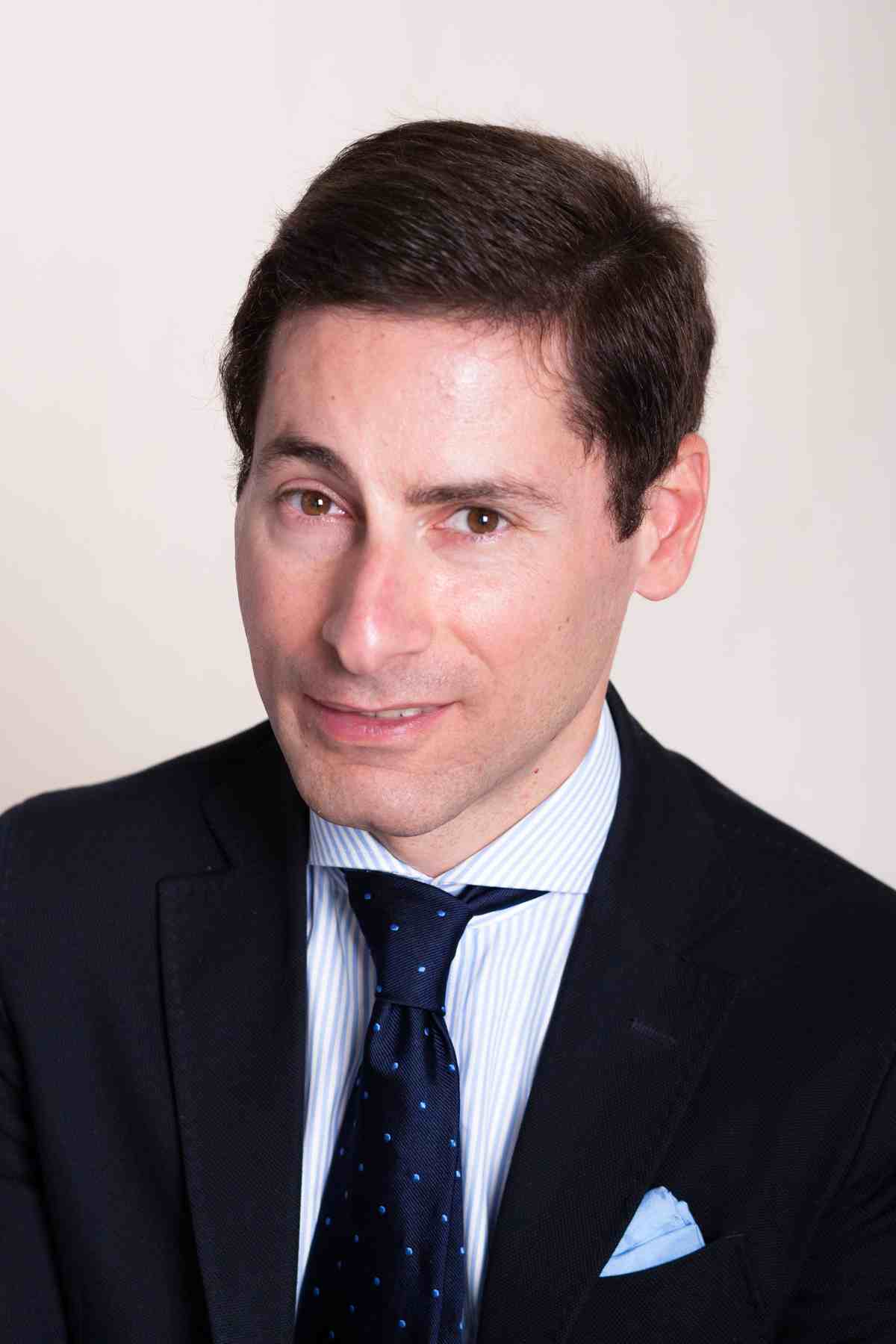}}]{Carlo Fischione} 
 is currently a Full Professor at KTH Royal Institute of Technology, Electrical Engineering and Computer Science, Stockholm, Sweden. He received the Ph.D. degree in Electrical and Information Engineering (3/3 years) in May 2005 from University of L'Aquila, Italy, and the Laurea degree in Electronic Engineering (Laurea, Summa cum Laude, 5/5 years) in April 2001 from the same University. He has held research positions at Massachusetts Institute of Technology, Cambridge, MA (2015, Visiting Professor); Harvard University, Cambridge, MA (2015, Associate); University of California at Berkeley, CA (2004-2005, Visiting Scholar, and 2007-2008, Research Associate). His research interests include optimization with applications to machine learning over networks, wireless sensor networks, networked control systems, and wireless networks. He received or co-received a number of awards, including the best paper award from the IEEE Transactions on Industrial Informatics (2007) and IEEE Transactions on Communications (2018).  He is Member of IEEE (the Institute of Electrical and Electronic Engineers), and Ordinary Member of DASP (the academy of history Deputazione Abruzzese di Storia Patria).
\end{IEEEbiography}

 \begin{IEEEbiography}[{\includegraphics[width=1in,height=1.25in,clip,keepaspectratio]{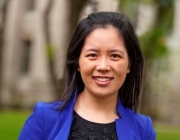}}]{Na Li} 
received her B.S. degree in mathematics and applied mathematics from Zhejiang University in China and her PhD degree in Control and Dynamical systems from the California Institute of Technology in 2013. She is an Associate Professor in the School of Engineering and Applied Sciences in Harvard University. She was a postdoctoral associate of the Laboratory for Information and Decision Systems at Massachusetts Institute of Technology.    

Her research lies in the design, analysis, optimization, and control of distributed network systems, with particular applications to cyber-physical network systems. She received NSF CAREER Award in 2016, AFOSR Young Investigator Award in 2017, ONR Young Investigator Award in 2019.
\end{IEEEbiography}

\end{document}